\documentclass{amsart}
\usepackage{graphicx}
\usepackage{amscd}
\usepackage{amsmath}
\usepackage{amsxtra}
\usepackage{amsfonts}
\usepackage{amssymb}
\newcommand{\volapukdef}{\newcommand{\volapuk}{\operatorname*{Const.}}}
\newtheorem{theorem}{Theorem}[section]
\newtheorem{corollary}[theorem]{Corollary}
\newtheorem{lemma}[theorem]{Lemma}
\newtheorem{proposition}[theorem]{Proposition}
\theoremstyle{definition}

\newtheorem{remark}[theorem]{Remark}

\theoremstyle{remark}
\newtheorem*{acknowledgements}{Acknowledgements}

\renewcommand{\theenumi}{\roman{enumi}}

\numberwithin{equation}{section}

\newcounter{linkequation}
\volapukdef
\makeatletter
\def\p@enumii{\relax}
\makeatother
\def\openone%{\hbox{\upshape \small1\kern-3.3pt\normalsize1}}
{\mathchoice
{\hbox{\upshape \small1\kern-3.3pt\normalsize1}}
{\hbox{\upshape \small1\kern-3.3pt\normalsize1}}
{\hbox{\upshape \tiny1\kern-2.3pt\SMALL1}}
{\hbox{\upshape \Tiny1\kern-2pt\tiny1}}}

\begin{document}
\title{Diagonalizing operators with reflection symmetry}
\author{Palle~E.~T.~ Jorgensen}
\address{Department of Mathematics\\
The University of Iowa\\
Iowa City, IA 52242-1419\\
U.S.A.}
\email{jorgen@math.uiowa.edu}
\thanks{Work supported in part by the National Science Foundation.}
\dedicatory{Dedicated to the memory of I.E. Segal}

\begin{abstract}
Let $U$ be an operator in a Hilbert space $\mathcal{H}_{0}$, and let
$\mathcal{K}\subset\mathcal{H}_{0}$ be a closed and invariant subspace.
Suppose there is a period-$2$ unitary operator $J$ in $\mathcal{H}_{0}$ such
that $JUJ=U^{\ast}$, and $PJP\geq0$, where $P$ denotes the projection of
$\mathcal{H}_{0}$ onto $\mathcal{K}$. We show that there is then a Hilbert
space $\mathcal{H}\left(  \mathcal{K}\right)  $, a contractive operator
$W\colon\mathcal{K}\rightarrow\mathcal{H}\left(  \mathcal{K}\right)  $, and a
selfadjoint operator $S=S\left(  U\right)  $ in $\mathcal{H}\left(
\mathcal{K}\right)  $ such that $W^{\ast}W=PJP$, $W$ has dense range, and
$SW=WUP$. Moreover, given $\left(  \mathcal{K},J\right)  $ with the stated
properties, the system $\left(  \mathcal{H}\left(  \mathcal{K}\right)
,W,S\right)  $ is unique up to unitary equivalence, and subject to the three
conditions in the conclusion. We also provide an operator-theoretic model of
this structure where $U|_{\mathcal{K}}$ is a pure shift of infinite
multiplicity, and where we show that $\ker\left(  W\right)  =0$. For that
case, we describe the spectrum of the selfadjoint operator $S\left(  U\right)
$ in terms of structural properties of $U$. In the model, $U$ will be realized
as a unitary scaling operator of the form%
\[
f\left(  x\right)  \longmapsto f\left(  cx\right)  ,\qquad c>1,
\]
and the spectrum of $S\left(  U_{c}\right)  $ is then computed in terms of the
given number $c$.
\end{abstract}\maketitle

\section{\label{Int}Introduction}

The paper is motivated by two problems one from mathematical physics, and the
other from the interface of integral transforms and interpolation theory. The
first problem is that of changing the spectrum of an operator, or a
one-parameter group of operators, with a view to getting a new spectrum with
physical desiderata (see, e.g., \cite{Seg98}), for example creating a mass
gap, and still preserving quasi-equivalence of the two underlying operator
systems. In the other problem we study how Hilbert space functional
completions change under the variation of a parameter in the integral kernel
of the transform in question. The motivating example here is derived from a
certain version of the Segal--Bargmann transform. For more detail on the
background and the applications alluded to in the Introduction, we refer to
the two previous joint papers \cite{JoOl98} and \cite{JoOl99}, as well as
\cite{Nee94} and \cite{Hal98}.

Let $U$ be an operator in a Hilbert space $\mathcal{H}_{0}$, and let $J$ be a
period-$2$ unitary operator in $\mathcal{H}_{0}$ such that%
\begin{equation}
JUJ=U^{\ast}. \label{eqInt.1}%
\end{equation}
We think of (\ref{eqInt.1}) as a reflection symmetry for the given operator
$U$. In this case, $U$ and its adjoint $U^{\ast}$ have the same spectrum, but,
of course, $U$ need not be selfadjoint. Nonetheless, we shall think of
(\ref{eqInt.1}) as a notion which generalizes selfadjointness. As an example,
let the Hilbert space $\mathcal{H}_{0}=L^{2}\left(  \mathbb{T}\right)  $,
\begin{equation}
\left(  Uf\right)  \left(  z\right)  =zf\left(  z\right)  ,\qquad f\in
L^{2}\left(  \mathbb{T}\right)  ,\;z\in\mathbb{T}, \label{eqInt.2ins}%
\end{equation}
and%
\begin{equation}
Jf\left(  z\right)  =f\left(  \bar{z}\right)  . \label{eqInt.3ins}%
\end{equation}
The space $L^{2}\left(  \mathbb{T}\right)  $ is from Haar measure on the
circle group $\mathbb{T}=\left\{  z\in\mathbb{C}\mathrel{;}\left|  z\right|
=1\right\}  $. It clear that (\ref{eqInt.1}) then holds. If $\mathcal{K}%
=H^{2}\left(  \mathbb{T}\right)  $ is the Hardy space of functions, $f\left(
z\right)  =\sum_{n=0}^{\infty}c_{n}z^{n}$, with $\left\|  f\right\|  ^{2}%
=\sum_{n=0}^{\infty}\left|  c_{n}\right|  ^{2}<\infty$, then we also have%
\begin{equation}
PJP\geq0 \label{eqInt.2bis}%
\end{equation}
where $P$ denotes the projection onto $H^{2}\left(  \mathbb{T}\right)  $. In
fact%
\begin{equation}
\left\langle f,Jf\right\rangle =\left|  c_{0}\right|  ^{2}, \label{eqInt.3bis}%
\end{equation}
where $\left\langle \,\cdot\,,\,\cdot\,\right\rangle $ denotes the inner
product in $L^{2}\left(  \mathbb{T}\right)  $. While our result applies to the
multiplicity-one shift, this is a degenerate situation, and the nontrivial
applications are for the case of infinite multiplicity.

There is in fact an infinite-multiplicity version of the above which we
proceed to describe. Let $0<s<1$ be given, and let $\mathcal{H}_{s}$ be the
Hilbert space whose norm $\left\|  f\right\|  _{s}$ is given by%
\begin{equation}
\left\|  f\right\|  _{s}^{2}=\int_{\mathbb{R}}\int_{\mathbb{R}}\overline
{f\left(  x\right)  }\,\left|  x-y\right|  ^{s-1}f\left(  y\right)  \,dx\,dy.
\label{eqInt.4}%
\end{equation}
Let $a\in\mathbb{R}_{+}$ be given, and set%
\begin{equation}
\left(  U\left(  a\right)  f\right)  \left(  x\right)  =a^{s+1}f\left(
a^{2}x\right)  . \label{eqInt.5}%
\end{equation}
It is clear that then $a\mapsto U\left(  a\right)  $ is a unitary
representation of the multiplicative group $\mathbb{R}_{+}$ acting on the
Hilbert space $\mathcal{H}_{s}$. It can be checked that $\left\|  f\right\|
_{s}$ in (\ref{eqInt.4}) is finite for all $f\in C_{c}\left(  \mathbb{R}%
\right)  $ ($=$ the space of compactly supported functions on the line). Now
let $\mathcal{K}$ ($=\mathcal{K}_{s}$) be the closure of $C_{c}\left(
-1,1\right)  $ in $\mathcal{H}_{s}$ relative to the norm $\left\|
\,\cdot\,\right\|  _{s}$ of (\ref{eqInt.4}). It is then immediate that
$U\left(  a\right)  $, for $a>1$, leaves $\mathcal{K}_{s}$ invariant, i.e., it
restricts to a semigroup of isometries $\left\{  U\left(  a\right)  \mathrel
{;}a>1\right\}  $ acting on $\mathcal{K}_{s}$. Setting%
\begin{equation}
\left(  Jf\right)  \left(  x\right)  =\left|  x\right|  ^{-s-1}f\left(
\frac{1}{x}\right)  ,\qquad x\in\mathbb{R}\setminus\left\{  0\right\}  ,
\label{eqInt.6}%
\end{equation}
we check that $J$ is then a period-$2$ unitary in $\mathcal{H}_{s}$, and that
\begin{equation}
JU\left(  a\right)  J=U\left(  a\right)  ^{\ast}=U\left(  a^{-1}\right)
\label{eqInt.7}%
\end{equation}
and%
\begin{equation}
\left\langle f,Jf\right\rangle _{\mathcal{H}_{s}}\geq0,\qquad\forall
\,f\in\mathcal{K}_{s}, \label{eqInt.8}%
\end{equation}
where $\left\langle \,\cdot\,,\,\cdot\,\right\rangle _{\mathcal{H}_{s}}$ is
the inner product
\begin{equation}
\left\langle f_{1},f_{2}\right\rangle _{\mathcal{H}_{s}}:=\int_{\mathbb{R}%
}\int_{\mathbb{R}}\overline{f_{1}\left(  x\right)  }\,\left|  x-y\right|
^{s-1}f_{2}\left(  y\right)  \,dx\,dy. \label{eqInt.9}%
\end{equation}
In fact, if $f\in C_{c}\left(  -1,1\right)  $, the expression in
(\ref{eqInt.8}) works out as the following reproducing kernel integral:%
\begin{equation}
\int_{-1}^{1}\int_{-1}^{1}\overline{f\left(  x\right)  }\left(  1-xy\right)
^{s-1}f\left(  y\right)  \,dx\,dy, \label{eqInt.10}%
\end{equation}
and we refer to \cite{JoOl98,JoOl99} for more details on this example.

As an application of our result, we will show that, if $a>1$, then $U\left(
a\right)  |_{\mathcal{K}_{s}}$ induces a selfadjoint operator $S\left(
a\right)  $ in a Hilbert space $\mathcal{H}\left(  \mathcal{K}_{s}\right)  $,
and there is a contraction $W\colon\mathcal{K}_{s}\rightarrow\mathcal{H}%
\left(  \mathcal{K}_{s}\right)  $, with%
\begin{equation}
\ker\left(  W\right)  =0, \label{eqInt.11}%
\end{equation}
such that%
\begin{equation}
W^{\ast}W=PJP, \label{eqInt.12}%
\end{equation}%
\begin{equation}
S\left(  a\right)  W=WU\left(  a\right)  P, \label{eqInt.13}%
\end{equation}
and%
\begin{equation}
\operatorname*{spectrum}\left(  S\left(  a\right)  \right)  =\left\{
a^{s-1-2n}\mathrel{;}n=0,1,2,\dots\right\}  . \label{eqInt.14}%
\end{equation}
What is important in this application is the property (\ref{eqInt.11}). So the
properties in this case for $W$ are $\left\|  W\right\|  \leq1$, $\ker\left(
W^{\ast}\right)  =\ker\left(  W\right)  =0$. While of course $U\left(
a\right)  |_{\mathcal{K}_{s}}$ and $S\left(  a\right)  $ cannot be unitarily
equivalent, then $W$ nonetheless defines a strong notion of equivalence
(quasi-equivalence) for the two semigroups $U\left(  a\right)  |_{\mathcal{K}%
_{s}}$ and $S\left(  a\right)  $, $a>1$, specified by the intertwining
property%
\begin{equation}
S\left(  a\right)  W=WU\left(  a\right)  P. \label{eqInt.15}%
\end{equation}
In particular, since both $W$ and $W^{\ast}$ have dense range in the
respective Hilbert spaces $\mathcal{K}$ and $\mathcal{H}\left(  \mathcal{K}%
\right)  $, it follows that the partial isometry part $L$ in the polar
decomposition $W=L\left(  W^{\ast}W\right)  ^{1/2}=L\left(  PJP\right)
^{1/2}$, is in fact a \emph{unitary }isomorphism of $\mathcal{K}$ onto
$\mathcal{H}\left(  \mathcal{K}\right)  $. The intertwining property for
$W^{\ast}W$ of the polar decomposition is%
\begin{equation}
\left(  W^{\ast}W\right)  UP=PU^{\ast}\left(  W^{\ast}W\right)  .
\label{eqInt.16}%
\end{equation}
But this cannot be iterated, so there is \emph{not} an analogous relation for
the factors $\left(  W^{\ast}W\right)  ^{1/2}$ and $L$. The properties of $W$
and $S$ in this example imply that $UP$ is in fact a pure shift (i.e., the
unitary part of the isometry $U|_{\mathcal{K}_{s}}$ of the Wold decomposition
is trivial, and moreover the backwards shift $PU^{\ast}$ has a cyclic vector.
The second conclusion is unique to this example, and follows from the fact
that $S=S\left(  a\right)  $ has simple spectrum.

\begin{proposition}
\label{ProInt.1}The isometry $UP$ is a pure shift.
\end{proposition}

\begin{proof}
The result may be read off from the following estimate:%
\begin{equation}
\left\|  PU^{\ast\,k}W^{\ast}\psi\right\|  =\left\|  W^{\ast}S\left(
a^{k}\right)  \psi\right\|  \leq\left\|  S\left(  a^{k}\right)  \psi\right\|
\leq a^{k\left(  s-1\right)  }\left\|  \psi\right\|  \underset{k\rightarrow
\infty}{\longrightarrow}0, \label{eqInt.17}%
\end{equation}
the estimate being valid for all $\psi\in\mathcal{H}\left(  \mathcal{K}%
\right)  $. Since $\ker\left(  W\right)  =0$, $W^{\ast}\mathcal{H}\left(
\mathcal{K}\right)  $ is dense in $\mathcal{K}$, so we have $\lim
_{k\rightarrow\infty}\left\|  PU^{\ast\,k}\varphi\right\|  =0$ for all
$\varphi\in\mathcal{K}$, and this last property is equivalent to
$U|_{\mathcal{K}_{s}}$ being a pure shift on $\mathcal{K}_{s}$.
\end{proof}

The restriction on $s$ remains $0<s<1$. It follows in fact from
\cite{JoOl98,JoOl99} that the multiplicity of this shift is $\infty$, i.e.,
that if $a>1$, the dimension of $\mathcal{K}_{s}\ominus U\left(  a\right)
\mathcal{K}_{s}$ is infinite.

The simplest case of a system $\left(  \mathcal{H}_{0},J\right)  $ with $J$ as
a reflection is that of $\mathcal{H}_{0}=\mathcal{H}\oplus\mathcal{H}$ and
$J=I\oplus\left(  -I\right)  $, i.e.,
\begin{equation}
J\left(  h_{1}\oplus h_{2}\right)  =h_{1}\oplus\left(  -h_{2}\right)  ,\qquad
h_{1},h_{2}\in\mathcal{H}. \label{eqInt.a}%
\end{equation}
In many applications of this, it will further be given that $\mathcal{H}$ is a
reproducing kernel Hilbert space in the sense of \cite{Aro50}. Suppose this is
the case, and that $Q\left(  \,\cdot\,,\,\cdot\,\right)  $ is the
corresponding reproducing kernel. We then have $\mathcal{H}$ realized as a
Hilbert space of $\mathbb{C}$-valued functions $h\left(  \,\cdot\,\right)  $
defined on some set $\Omega$, and $Q$ is a function on $\Omega\times\Omega$
such that $Q\left(  z,\,\cdot\,\right)  \in\mathcal{H}$ for all $z\in\Omega$,
and
\begin{equation}
\left\langle Q\left(  z,\,\cdot\,\right)  ,h\right\rangle =h\left(  z\right)
\text{\qquad for all }h\in\mathcal{H}. \label{eqInt.b}%
\end{equation}
In this case, we will use $Q$ in identifying a class of subspaces
$\mathcal{K}\subset\mathcal{H}\oplus\mathcal{H}$ such that%
\begin{equation}
\left\langle k,Jk\right\rangle \geq0\text{\qquad for all }k\in\mathcal{K}.
\label{eqInt.c}%
\end{equation}

We now describe such a class of spaces $\mathcal{K}$. Let $D:=\left\{
z\in\mathbb{C}\mathrel{;}\left|  z\right|  <1\right\}  $. It will be stated in
an abstract setting, and the applications to interpolation theory will be
given in Section \ref{Han} below.

\begin{proposition}
\label{ProInt.3}Let $\mathcal{H}$ be a reproducing kernel Hilbert space
corresponding to a kernel function%
\begin{equation}
Q\colon\Omega\times\Omega\longrightarrow\mathbb{C}, \label{eqInt.d}%
\end{equation}
and let $\Omega_{0}\subset\Omega$ be a subset. Let a function
\begin{equation}
\varphi\colon\Omega_{0}\longrightarrow\bar{D} \label{eqInt.e}%
\end{equation}
be given, and let $\mathcal{K}_{\varphi}\subset\mathcal{H}\oplus\mathcal{H}$
be defined as the closed span of%
\begin{equation}
\left\{
\begin{pmatrix}
Q\left(  z,\,\cdot\,\right) \\
\varphi\left(  z\right)  Q\left(  z,\,\cdot\,\right)
\end{pmatrix}
\mathrel{;}z\in\Omega_{0}\right\}  \subset%
\begin{pmatrix}
\mathcal{H}\\
\mathcal{H}%
\end{pmatrix}
^{\oplus}. \label{eqInt.f}%
\end{equation}

\begin{enumerate}
\item \label{ProInt.3(1)}Then \textup{(\ref{eqInt.c})} holds for $J=\left(
\begin{smallmatrix}
I & 0\\
0 & -I
\end{smallmatrix}
\right)  $ if and only if
\[
\left(  z_{1},z_{2}\right)  \longmapsto\left(  1-\overline{\varphi\left(
z_{1}\right)  }\varphi\left(  z_{2}\right)  \right)  Q\left(  z_{1}%
,z_{2}\right)
\]
is positive definite on $\Omega_{0}$.

\item \label{ProInt.3(2)}If instead $\varphi\colon\Omega_{0}\rightarrow
\mathbb{C}$, and $J=\left(
\begin{smallmatrix}
0 & I\\
I & 0
\end{smallmatrix}
\right)  $, then \textup{(\ref{eqInt.c})} holds if and only if
\[
\left(  z_{1},z_{2}\right)  \longmapsto\left(  \overline{\varphi\left(
z_{1}\right)  }+\varphi\left(  z_{2}\right)  \right)  Q\left(  z_{1}%
,z_{2}\right)
\]
is positive definite on $\Omega_{0}$.
\end{enumerate}
\end{proposition}

\begin{proof}
The result follows from a substitution of the vectors in (\ref{eqInt.f}) into
the positivity requirement (\ref{eqInt.c}), and computing out the answer for
the two cases of reflection $J$, i.e., $J=\left(
\begin{smallmatrix}
I & 0\\
0 & -I
\end{smallmatrix}
\right)  $ and $J=\left(
\begin{smallmatrix}
0 & I\\
I & 0
\end{smallmatrix}
\right)  $. We refer to Section \ref{Han} for more details, and additional
comments on applications to interpolation theory.
\end{proof}

\section{\label{Pur}Pure isometries}

It is well known that pure isometries (alias shifts) of infinite multiplicity
play a role in the harmonic analysis of wavelets, see \cite{BrJo97b}, and in
the Lax--Phillips version of scattering theory for the wave equation
\cite{LaPh89}. Let $V$ be a shift in a Hilbert space $\mathcal{K}$, and let%
\begin{equation}
\mathcal{L}:=\mathcal{K}\ominus V\mathcal{K}; \label{eqPur.1}%
\end{equation}
then%
\begin{equation}
\mathcal{K}=\sideset{}{^{\smash{\oplus}}}{\sum}\limits_{n=0}^{\infty}%
V^{n}\mathcal{L} \label{eqPur.2}%
\end{equation}
as a direct sum. But for every nonzero $l\in\mathcal{L}$, and $z\in
D:=\left\{  z\in\mathbb{C}\mathrel{;}\left|  z\right|  <1\right\}  $, the
vector%
\begin{equation}
f=f\left(  l,z\right)  :=l\oplus zVl\oplus z^{2}V^{2}l\oplus\cdots
\label{eqPur.3}%
\end{equation}
is an eigenvector of $V^{\ast}$, i.e.,%
\begin{equation}
V^{\ast}f=zf, \label{eqPur.4}%
\end{equation}
and $\left\|  f\right\|  ^{2}=\left(  1-\left|  z\right|  ^{2}\right)
^{-1}\left\|  l\right\|  ^{2}$. In fact, as $l$ varies over $\mathcal{L}%
\setminus\left\{  0\right\}  $, the vectors
\begin{equation}
\left\{  f\left(  l,z^{n}\right)  \mathrel{;}n=1,2,\dots\right\}
\label{eqPur.5}%
\end{equation}
span a dense subspace in $\mathcal{K}$. This is true for every $z\in D$ fixed;
so it is clear from this that there is a variety of ways of creating
selfadjoint, and normal, realizations of a given $V$, i.e., solutions to the
problem%
\begin{equation}
WV=NW. \label{eqPur.6}%
\end{equation}
Specifically, there is a Hilbert space $\mathcal{H}\left(  \mathcal{K}\right)
$, a bounded operator $W\colon\mathcal{K}\rightarrow\mathcal{H}\left(
\mathcal{K}\right)  $, and a normal operator $N$ in $\mathcal{H}\left(
\mathcal{K}\right)  $ such that (\ref{eqPur.6}) holds. This problem has been
studied recently by Feldman \cite{Fel99}, and Agler et al.\ \cite{AgMc98}, but
it is a different focus from ours. The reflection $J$ plays a crucial role in
our approach. It also makes our setting considerably more restrictive and it
allows us to get solutions to the diagonalization problem which are
\emph{unique up to unitary equivalence.} More importantly, it gives an answer
to a reflection problem from mathematical physics which we proceed to describe.

The approach (\ref{eqPur.4}) for $V^{\ast}$ works for a wider class of
operators than the backwards shift, namely the operators in the Cowen--Douglas
classes, see \cite{CoDo78}, but we have not yet checked which of the
Cowen--Douglas operators admit reflection symmetry.

Our next result will be stated for general bounded operators $U$ which have
reflection symmetry, and the symmetry is given in terms of a period-$2$
unitary $J$ and a subspace $\mathcal{K}$ which is invariant under $U$. From
this we will then arrive at a selfadjoint realization $S$ of $U$, and when
$\left(  \mathcal{K},J\right)  $ is given, we will show that $S$ is determined
uniquely up to unitary equivalence. The result is interesting even if $U$ is
given at the outset to be unitary. In fact in an application from quantum
field theory, $U$ will be rather a unitary one-parameter group $\left\{
U\left(  t\right)  \right\}  _{t\in\mathbb{R}}$ of operators acting on a
Hilbert space $\mathcal{H}_{0}$, and $\mathcal{K}$ will be a subspace in
$\mathcal{H}_{0}$ which is invariant under $U\left(  t\right)  $ for $t\geq0$.
By Stone's theorem \cite{Var85}, there is a selfadjoint Hamiltonian operator
$H$ (generally unbounded) in $\mathcal{H}_{0}$ such that%
\begin{equation}
U\left(  t\right)  =e^{-itH},\qquad t\in\mathbb{R}. \label{eqPur.7}%
\end{equation}
In this application, we will have%
\begin{equation}
JU\left(  t\right)  J=U\left(  -t\right)  ,\qquad t\in\mathbb{R},
\label{eqPur.8}%
\end{equation}
and $J$ is referred to as ``time-reversal'' or ``time-reflection''. The
initial Hamiltonian might not have the right ``physical'' spectrum; for
example, the spectrum of $H$ might be all of $\mathbb{R}$, and what is desired
would be a spectrum which is contained in $\mathbb{R}_{+}$ with a positive gap
between $0$ and the bottom of the ``physical'' spectrum. We will show that
this can be achieved; in fact we will describe a selfadjoint realization
$S=S\left(  U\right)  $ in the form of a semigroup%
\begin{equation}
S\left(  t\right)  =e^{-t\hat{H}} \label{eqPur.9}%
\end{equation}
where $\hat{H}$ is a selfadjoint operator in the new Hilbert space
$\mathcal{H}\left(  \mathcal{K}\right)  $, and the spectrum of $\hat{H}$ will
be ``physical'' in that it will be positive and there will be a ``mass gap'',
i.e., a positive gap between $0$ and the lower bound for
$\operatorname*{spectrum}\left(  \hat{H}\right)  $. But the key to passing
from $H$ to $\hat{H}$ will be the given $\left(  \mathcal{K},J\right)  $ when
$\mathcal{K}\subset\mathcal{H}_{0}$ is assumed invariant under $U\left(
t\right)  $, $t\geq0$, and $J$ is a time-reflection, i.e., $J$ and $\left\{
U\left(  t\right)  \right\}  $ will satisfy (\ref{eqPur.8}). As we noted, the
construction $H\rightsquigarrow\hat{H}$ with $\hat{H}$ having a mass-gap will
show, after the fact, that the initial semigroup of isometries $U\left(
t\right)  |_{\mathcal{K}}$, $t\geq0$, will necessarily be a pure shift (and of
infinite multiplicity). By this we mean that there is a unitary isomorphism
between $\mathcal{H}_{0}$ and $L^{2}\left(  \mathbb{R},\mathcal{M}\right)  $
for some infinite-dimensional Hilbert space $\mathcal{M}$ which intertwines
$\left\{  U\left(  t\right)  \right\}  _{t\in\mathbb{R}}$ with translation on
$L^{2}\left(  \mathbb{R},\mathcal{M}\right)  $. Specifically, there is a
unitary isomorphism%
\begin{equation}
Y\colon\mathcal{H}_{0}\longrightarrow L^{2}\left(  \mathbb{R},\mathcal{M}%
\right)  \text{, onto,} \label{eqPur.10}%
\end{equation}
such that%
\begin{equation}
YU\left(  t\right)  Y^{-1}f\left(  x\right)  =f\left(  x-t\right)  ,\qquad
f\in L^{2}\left(  \mathbb{R},\mathcal{M}\right)  ,\;t\in\mathbb{R},
\label{eqPur.11}%
\end{equation}
with the further property that%
\begin{equation}
Y\left(  \mathcal{K}\right)  =L^{2}\left(  \mathbb{R}_{+},\mathcal{M}\right)
, \label{eqPur.12}%
\end{equation}
i.e., the functions in $L^{2}\left(  \mathbb{R},\mathcal{M}\right)  $ which
are supported in the positive half line.

\section{\label{Ref}Reflection symmetry}

The following result provides the axiomatic setup for reflection symmetry in
the form described above. With the given symmetry axioms, it provides the step
$U\mapsto S\left(  U\right)  $ from a general operator $U$ with symmetry to
its selfadjoint version $S\left(  U\right)  $, and we show that $S\left(
U\right)  $ is unique up to unitary equivalence. The data that emerges is
$\left(  \mathcal{H}\left(  \mathcal{K}\right)  ,W,S\right)  $, where%
\begin{equation}
SW=WUP. \label{eqRef.1}%
\end{equation}
Here $P$ denotes the projection onto the subspace $\mathcal{K}$ which both is
invariant for $U$ and satisfies reflection positivity relative to the
period-$2$ unitary $J$ (i.e., the reflection). But in the general setting, the
axioms allow $W\colon\mathcal{K}\rightarrow\mathcal{H}\left(  \mathcal{K}%
\right)  $ to have nonzero kernel, and this represents some degree of
non-uniqueness: for example, $W$ may be a ``small'' (rank-one, say)
projection, and $S$ might be zero. Hence we shall focus on the setting when
$\ker\left(  W\right)  =0$, and we will say then that the two operators
$U|_{\mathcal{K}}$ and $S$ are \emph{quasi-equivalent.} While the intertwining
operator $W$ is $1$--$1$ with dense range, its inverse $W^{-1}$ will be
unbounded.\setcounter{linkequation}{\value{equation}}

\begin{theorem}
\label{ThmRef.1}\setcounter{equation}{0} \renewcommand{\theequation}%
{\roman{equation}} \renewcommand{\theenumi}{\alph{enumi}}\raggedright Let $U$
be a bounded operator in a Hilbert space $\mathcal{H}_{0}$. Let $\mathcal{K}%
\subset\mathcal{H}_{0}$ be an invariant subspace, and let $P$ denote the
projection of $\mathcal{H}_{0}$ onto $\mathcal{K}$. Let $J$ be a period-$2$
unitary operator in $\mathcal{H}_{0}$ which satisfies%
\begin{equation}
JUJ=U^{\ast} \label{eqThmRef.1(1)}%
\end{equation}
and%
\begin{equation}
PJP\geq0. \label{eqThmRef.1(2)}%
\end{equation}

\begin{enumerate}
\item \label{ThmRef.1(1)}Then there is a Hilbert space $\mathcal{H}\left(
\mathcal{K}\right)  $ and a contractive operator%
\[
W\colon\mathcal{K}\longrightarrow\mathcal{H}\left(  \mathcal{K}\right)
\]
with dense range, and a bounded selfadjoint operator $S=S\left(  U\right)  $
in $\mathcal{H}\left(  \mathcal{K}\right)  $ such that%
\begin{equation}
SW=WUP, \label{eqThmRef.1(3)}%
\end{equation}%
\begin{equation}
W^{\ast}W=PJP, \label{eqThmRef.1(4)}%
\end{equation}
and%
\begin{equation}
\left\|  S\left(  U\right)  \right\|  \leq\left(  \operatorname*{sp}\left(
U^{2}\right)  \right)  ^{\frac{1}{2}}, \label{eqThmRef.1(5)}%
\end{equation}
where $\operatorname*{sp}\left(  U^{2}\right)  $ denotes the spectral radius
of $U^{2}$.

\item \label{ThmRef.1(2)}Given \textup{(\ref{eqThmRef.1(1)}%
)--(\ref{eqThmRef.1(2)})} the data $\left(  \mathcal{H}\left(  \mathcal{K}%
\right)  ,W,S\right)  $ is unique up to unitary equivalence subject to the
axioms \textup{(\ref{eqThmRef.1(3)})--(\ref{eqThmRef.1(4)}).} Specifically,
suppose $\left(  \mathcal{H}_{i}\left(  \mathcal{K}\right)  ,W_{i}%
,S_{i}\right)  $, $i=1,2$, are two systems which both solve the extension
problem, i.e., are extensions satisfying \textup{(\ref{eqThmRef.1(3)}%
)--(\ref{eqThmRef.1(4)}).} Then there is a unitary isomorphism $T\colon
\mathcal{H}_{1}\left(  \mathcal{K}\right)  \rightarrow\mathcal{H}_{2}\left(
\mathcal{K}\right)  $ of $\mathcal{H}_{1}\left(  \mathcal{K}\right)  $ onto
$\mathcal{H}_{2}\left(  \mathcal{K}\right)  $ which satisfies
\begin{equation}
TW_{1}=W_{2} \label{eqThmRef.1(6)}%
\end{equation}
and%
\begin{equation}
TS_{1}=S_{2}T. \label{eqThmRef.1(7)}%
\end{equation}

\item \label{ThmRef.1(3)}There are operators $U$, with reflection symmetry,
such that $W$ from $\left(  \mathcal{H}\left(  \mathcal{K}\right)
,W,S\right)  $ has%
\begin{equation}
\ker\left(  W\right)  =0. \label{eqThmRef.1(8)}%
\end{equation}
\end{enumerate}
\end{theorem}

\setcounter{equation}{\value{linkequation}}

\begin{proof}
The proof is rather long and will be broken up into its three parts
(\ref{ThmRef.1(1)}), (\ref{ThmRef.1(2)}), and (\ref{ThmRef.1(3)}). Part
(\ref{ThmRef.1(1)}) asserts the existence of a selfadjoint realization of the
given operator $U$, while part (\ref{ThmRef.1(2)}) is uniqueness up to unitary
equivalence. Part (\ref{ThmRef.1(3)}) is an explicit construction which takes
place in a certain reproducing kernel Hilbert space.

The following observation gives a more concrete understanding of axiom
(\ref{eqThmRef.1(2)}) in part (\ref{ThmRef.1(1)}) of Theorem \ref{ThmRef.1}.
Let $J$ be a period-$2$ unitary operator in a Hilbert space $\mathcal{H}_{0}$,
and let $\mathcal{H}_{\pm}$ be the respective eigenspaces corresponding to
eigenvalues $\pm1$ of $J$. If $P_{+}$ is the projection onto $\mathcal{H}_{+}%
$, then $J=2P_{+}-I$.\renewcommand{\qed}{}
\end{proof}

\begin{lemma}
\label{LemRef.2}A closed subspace $\mathcal{K}\subset\mathcal{H}_{0}$
satisfies \textup{(\ref{eqThmRef.1(2)})} if and only if $\mathcal{K}$ is the
graph of a contractive operator $\Lambda$ from $\mathcal{H}_{+}$ to
$\mathcal{H}_{-}$. By this we mean that $\Lambda$ is defined on a closed
subspace $\mathcal{P}\subset\mathcal{H}_{+}$ and $\Lambda$ maps $\mathcal{P}$
contractively into $\mathcal{H}_{-}$. Hence $\mathcal{K}\simeq\left\{  \left(
p,\Lambda p\right)  \mathrel{;}p\in\mathcal{P}\right\}  $, or we will write
simply $\mathcal{K}=G\left(  \Lambda\right)  $ and $\mathcal{P}=D\left(
\Lambda\right)  $ where $G$ and $D$ are used for graph and domain, respectively.
\end{lemma}

\begin{proof}
The main idea in the proof is in \cite{Phil}, but we include a sketch. This
will also give us a chance for introducing some terminology which will be
needed later anyway. Suppose $\mathcal{K}\subset\mathcal{H}_{0}$ is a closed
subspace which satisfies (\ref{eqThmRef.1(2)}). For $k\in\mathcal{K}$ we have
$k=P_{+}k+P_{-}k$, where $P_{-}:=I-P_{+}$ and $J=P_{+}-P_{-}$. But
$\left\langle k,Jk\right\rangle =\left\|  P_{+}k\right\|  ^{2}-\left\|
P_{-}k\right\|  ^{2}$ for all $k\in\mathcal{K}$ by (\ref{eqThmRef.1(2)}), and
if we define $\Lambda P_{+}k:=P_{-}k$, then $\Lambda$ is well-defined and
contractive from $\mathcal{P}=P_{+}\mathcal{K}$ to $P_{-}\mathcal{K}$. The
reasoning shows that the converse argument is also valid, so the lemma follows
except for the assertion that $\mathcal{P}:=P_{+}\mathcal{K}$ must be
automatically closed. Let $k_{n}$ be a sequence of vectors in $\mathcal{K}$
such that $P_{+}k_{n}\rightarrow h_{+}\in\mathcal{H}_{+}$. Then by
(\ref{eqThmRef.1(2)}),%
\[
\left\|  P_{-}\left(  k_{n}-k_{m}\right)  \right\|  \leq\left\|  P_{+}\left(
k_{n}-k_{m}\right)  \right\|  \longrightarrow0\text{\qquad as }%
n,m\longrightarrow\infty.
\]
So the limit $\lim_{n\rightarrow\infty}P_{-}k_{n}=h_{-}$ exists in
$\mathcal{H}_{-}$, and
\[
k_{n}=P_{+}k_{n}+P_{-}k_{n}\longrightarrow h_{+}+h_{-}.
\]
Since $\mathcal{K}$ is assumed closed in $\mathcal{H}_{0}$, we get
$h_{+}+h_{-}\in\mathcal{K}$, and $h_{+}=P_{+}\left(  h_{+}+h_{-}\right)
=\lim_{n\rightarrow\infty}P_{+}k_{n}$. This shows that $P_{+}\mathcal{K}$ is
closed, and the proof is completed.
\end{proof}

\begin{proof}
[Proof of Theorem \textup{\ref{ThmRef.1}} continued](\ref{ThmRef.1(1)}) Let
the operator $U$ be given as in the statement of the theorem. Let
$\mathcal{K}\subset\mathcal{H}_{0}$ be the invariant subspace with projection
$P$, and let $J$ be the reflection. It is assumed to satisfy
(\ref{eqThmRef.1(1)})--(\ref{eqThmRef.1(2)}). In view of (\ref{eqThmRef.1(2)}%
), we have
\begin{equation}
\left\langle k,Jk\right\rangle \geq0\text{\qquad for all }k\in\mathcal{K},
\label{eqRef.2}%
\end{equation}
where $\left\langle \,\cdot\,,\,\cdot\,\right\rangle $ denotes the given inner
product from $\mathcal{H}_{0}$. (Note that $\mathcal{K}$ is not invariant
under $J$, so the vector $Jk$ is typically not in $\mathcal{K}$ if $k$ is.)
Applying the Cauchy--Schwarz inequality, we get%
\begin{equation}
\left|  \left\langle k_{1},Jk_{2}\right\rangle \right|  ^{2}\leq\left\langle
k_{1},Jk_{1}\right\rangle \left\langle k_{2},Jk_{2}\right\rangle \text{\qquad
for all }k_{1},k_{2}\in\mathcal{K}. \label{eqRef.3}%
\end{equation}
The idea is to get a new Hilbert space $\mathcal{H}\left(  \mathcal{K}\right)
$ from the form $\left\langle k_{1},Jk_{2}\right\rangle $, i.e., that this
form should be the new inner product. So we must form the quotient space
$\mathcal{K}/\mathcal{N}$ where%
\begin{equation}
\mathcal{N}=\left\{  k\in\mathcal{K}\mathrel{;}\left\langle k,Jk\right\rangle
=0\right\}  . \label{eqRef.4}%
\end{equation}
In view of (\ref{eqRef.3}), we get%
\begin{equation}
\mathcal{N}=\left\{  k_{0}\in\mathcal{K}\mathrel{;}\left\langle k_{0}%
,Jk\right\rangle =0\text{ for all }k\in\mathcal{K}\right\}  . \label{eqRef.5}%
\end{equation}
Since%
\begin{equation}
\left\langle k_{1},JUk_{2}\right\rangle =\left\langle k_{1},U^{\ast}%
Jk_{2}\right\rangle =\left\langle Uk_{1},Jk_{2}\right\rangle \text{\qquad for
all }k_{1},k_{2}\in\mathcal{K}, \label{eqRef.6}%
\end{equation}
we conclude that $U$ passes to the quotient $\mathcal{K}/\mathcal{N}$ and
defines there a symmetric operator. When $\mathcal{K}/\mathcal{N}$ is
completed in the new norm $\left\|  \,\cdot\,\right\|  _{J}$,%
\begin{equation}
\left\|  k\right\|  _{J}^{2}:=\left\langle k,Jk\right\rangle , \label{eqRef.7}%
\end{equation}
the induced operator becomes selfadjoint in this Hilbert space%
\begin{equation}
\mathcal{H}\left(  \mathcal{K}\right)  :=\left(  \mathcal{K}/\mathcal{N}%
\right)  \sptilde. \label{eqRef.8}%
\end{equation}
The induced operator will be denoted $S=S\left(  U\right)  $, and we will now
show that it satisfies conditions (\ref{eqThmRef.1(3)})--(\ref{eqThmRef.1(5)}%
), starting with (\ref{eqThmRef.1(5)}), i.e., showing first that $S\left(
U\right)  $ is a bounded operator in the Hilbert space $\mathcal{H}\left(
\mathcal{K}\right)  $. The argument for boundedness is essentially in
\cite{JoOl98}, but we include it here for the convenience of the reader.

Let $k\in\mathcal{K}$, and use recursion on (\ref{eqRef.3}) as follows:%
\begin{align*}
\left\|  Uk\right\|  _{J}^{2}  &  =\left\langle Uk,JUk\right\rangle
=\left\langle Uk,U^{\ast}Jk\right\rangle =\left\langle U^{2}k,Jk\right\rangle
\\
&  \leq\left\langle U^{2}k,JU^{2}k\right\rangle ^{\frac{1}{2}}\left\langle
k,Jk\right\rangle ^{\frac{1}{2}}\\
&  \leq\left\langle U^{4}k,JU^{4}k\right\rangle ^{\frac{1}{4}}\left\langle
k,Jk\right\rangle ^{\frac{1}{2}+\frac{1}{4}}\\
&  \leq\vphantom{\left\langle U^{4}k,JU^{4}k\right\rangle^{\frac{1}{4}}%
\left\langle k,Jk\right\rangle^{\frac{1}{2}+\frac{1}{4}}}\cdots\\
&  \leq\left\langle U^{2^{n}}k,JU^{2^{n}}k\right\rangle ^{\frac{1}%
{2^{n\mathstrut}}}\cdot\left\langle k,Jk\right\rangle ^{\frac{1}%
{2^{\mathstrut}}+\frac{1}{4^{\mathstrut}}+\dots+\frac{1}{2^{n\mathstrut}}}\\
&  \leq\left\langle U^{2^{n+1}}k,Jk\right\rangle ^{\frac{1}{2^{n\mathstrut}}%
}\cdot\left\|  k\right\|  _{J}^{2}\\
&  \leq\left\|  U^{2^{n+1}}\right\|  ^{\frac{1}{2^{n\mathstrut}}}\cdot\left\|
k\right\|  ^{\frac{1}{2^{n-1\mathstrut}}}\cdot\left\|  k\right\|  _{J}^{2}.
\end{align*}
We have $\lim\limits_{n\rightarrow\infty}\left\|  U^{2^{n+1}}\right\|
^{\frac{1}{2^{n\mathstrut}}}=\operatorname*{sp}\left(  U^{2}\right)  =$ the
spectral radius, and $\lim\limits_{n\rightarrow\infty}\left\|  k\right\|
^{\frac{1}{2^{n-1\mathstrut}}}=1$ if $k\neq0$. We have therefore proved the
estimate%
\[
\left\|  Uk\right\|  _{J}\leq\left(  \operatorname*{sp}\left(  U^{2}\right)
\right)  ^{\frac{1}{2}}\left\|  k\right\|  _{J}%
\]
for $k\in\mathcal{K}$, and it follows that the induced operator $S=S\left(
U\right)  $ on $\mathcal{H}\left(  \mathcal{K}\right)  =\left(  \mathcal{K}%
/\mathcal{N}\right)  \sptilde$ satisfies (\ref{eqThmRef.1(5)}), as claimed.
Since we already showed that $S$ is selfadjoint, we conclude that $S$ has
bounded spectrum inside the interval%
\begin{equation}
\left[  -\left(  \operatorname*{sp}\left(  U^{2}\right)  \right)  ^{\frac
{1}{2}},\left(  \operatorname*{sp}\left(  U^{2}\right)  \right)  ^{\frac{1}%
{2}}\right]  \subset\mathbb{R}. \label{eqRef.9}%
\end{equation}
If $U$ on $\mathcal{H}_{0}$ is unitary, this is the interval $\left[
-1,1\right]  $. If $U=U\left(  t\right)  $, $t\in\mathbb{R}$, is a group of
operators, then $S=S\left(  t\right)  $, $t\geq0$, is a semigroup of
selfadjoint operators, and so%
\begin{equation}
S\left(  t\right)  =S\left(  \frac{t}{2}\right)  ^{2}\geq0 \label{eqRef.10}%
\end{equation}
for all $t\geq0$, and the spectrum of $S\left(  t\right)  $ is therefore
positive in that case, and we get the representation%
\begin{equation}
S\left(  t\right)  =e^{-t\hat{H}},\qquad t\geq0, \label{eqRef.11}%
\end{equation}
for some (generally unbounded) selfadjoint operator $\hat{H}$ in
$\mathcal{H}\left(  \mathcal{K}\right)  $.
\end{proof}

\begin{proof}
[Proof of part \textup{(\ref{ThmRef.1(2)})}]For a given operator $U$ which has
a pair $\left(  \mathcal{K},J\right)  $ defining a reflection symmetry, we
showed in (\ref{ThmRef.1(1)}) that there is a system $\left(  \mathcal{H}%
\left(  \mathcal{K}\right)  ,W,S\right)  $ with a selfadjoint operator $S$ in
$\mathcal{H}\left(  \mathcal{K}\right)  $, and an intertwining operator $W$,
which satisfy (\ref{eqThmRef.1(3)})--(\ref{eqThmRef.1(5)}) in the statement of
the theorem. We now prove that this system is unique up to unitary
equivalence. So suppose there are two systems $\left(  \mathcal{H}_{i}\left(
\mathcal{K}\right)  ,W_{i},S_{i}\right)  $, $i=1,2$, both satisfying
(\ref{eqThmRef.1(3)})--(\ref{eqThmRef.1(4)}) and with the two ``extension''
operators $S_{1}$ and $S_{2}$ both selfadjoint and bounded. We will now show
that there is then a unitary isomorphism $T\colon\mathcal{H}_{1}\left(
\mathcal{K}\right)  \rightarrow\mathcal{H}_{2}\left(  \mathcal{K}\right)  $
which defines the equivalence, i.e., which satisfies (\ref{eqThmRef.1(6)}) and
(\ref{eqThmRef.1(7)}) in the theorem. We will make (\ref{eqThmRef.1(6)}) into
a definition, setting%
\begin{equation}
TW_{1}k=W_{2}k, \label{eqRef.12}%
\end{equation}
for $k\in\mathcal{K}$. Since both $W_{1}$ and $W_{2}$ satisfy
(\ref{eqThmRef.1(4)}), we conclude that%
\[
\left\|  W_{1}k\right\|  _{J}=0\iff k\in\mathcal{N}\iff\left\|  W_{2}%
k\right\|  _{J}=0,
\]
or, stated equivalently,%
\[
\ker\left(  W_{i}\right)  =\mathcal{N}\text{\qquad for }i=1,2,
\]
where $\mathcal{N}$ is defined in (\ref{eqRef.4}). Hence, formula
(\ref{eqRef.12}) makes a good definition of a linear operator $T$ mapping a
dense subspace in $\mathcal{H}_{1}\left(  \mathcal{K}\right)  $ into one in
$\mathcal{H}_{2}\left(  \mathcal{K}\right)  $. But property
(\ref{eqThmRef.1(4)}) for $W_{1}$ and $W_{2}$ implies that $T$ is also
isometric, indeed%
\[
\left\|  TW_{1}k\right\|  _{J}^{2}=\left\|  W_{2}k\right\|  _{J}%
^{2}=\left\langle k,Jk\right\rangle =\left\|  W_{1}k\right\|  _{J}^{2}.
\]
Hence $T$ is a unitary isomorphism of $\mathcal{H}_{1}\left(  \mathcal{K}%
\right)  $ onto $\mathcal{H}_{2}\left(  \mathcal{K}\right)  $. Using now
(\ref{eqThmRef.1(3)}) for the two systems, we get%
\[
\left(  TS_{1}\right)  W_{1}k=TW_{1}Uk=W_{2}Uk=S_{2}W_{2}k=\left(
S_{2}T\right)  W_{1}k\text{\qquad for all }k\in\mathcal{K}.
\]
Since $W_{1}$ has dense range, we get the desired intertwining property
(\ref{eqThmRef.1(7)}) as claimed in the theorem.
\end{proof}

\begin{proof}
[Proof of part \textup{(\ref{ThmRef.1(3)})}]The assertion in part
(\ref{ThmRef.1(3)}) is that there are examples where the induction
$U\rightsquigarrow S\left(  U\right)  $ has intertwining operator $W$ with
zero kernel, or equivalently, $\mathcal{N}=\left\{  0\right\}  $. We already
mentioned this in (\ref{eqInt.4})--(\ref{eqInt.6}) of Section \ref{Int}, and
in fact this is a one-parameter semigroup of isometries $U\left(  a\right)
P_{\mathcal{K}_{s}}$, $a>1$. In fact, it arises as the restriction to an
invariant subspace of a unitary one-parameter group. It is a representation
$U\left(  a\right)  $, $a\in\mathbb{R}_{+}$, of the multiplicative group
$\mathbb{R}_{+}$, or equivalently, via $a=e^{t}$, a representation of the
additive group $\mathbb{R}$. We get as a corollary of (\ref{ThmRef.1(3)}) that
$\left\{  U_{s}\left(  e^{t}\right)  \right\}  _{t\in\mathbb{R}}$ is
equivalent to the group of translations on $L^{2}\left(  \mathbb{R}%
,\mathcal{M}\right)  $ for some infinite-dimensional Hilbert space
$\mathcal{M}$ as described in (\ref{eqPur.10})--(\ref{eqPur.12}) in the
conclusion of Section \ref{Pur} above.

Now recall the Hilbert space $\mathcal{H}_{s}$ and its subspace $\mathcal{K}%
_{s}$ from Section \ref{Int}. When $0<s<1$, $\mathcal{H}_{s}$ is defined by
the norm $\left\|  \,\cdot\,\right\|  _{s}$ from (\ref{eqInt.4}) and the
subspace $\mathcal{K}_{s}$ is the completion of $C_{c}\left(  -1,1\right)  $
in the $\left\|  \,\cdot\,\right\|  _{s}$-norm. We may pick some $a>1$, and
consider the isometry $U_{s}\left(  a\right)  |_{\mathcal{K}_{s}}$ of
$\mathcal{K}_{s}$. From (\ref{eqInt.6}) we see that $J$ also depends on $s$.
The new inner product is
\begin{equation}
\left\langle k_{1},k_{2}\right\rangle _{J}:=\left\langle k_{1},Jk_{2}%
\right\rangle _{\mathcal{H}_{s}}\label{eqRef.13}%
\end{equation}
(defined for $k_{1},k_{2}\in\mathcal{K}_{s}$), and depends on $s$ as well. It
is worked out explicitly in (\ref{eqInt.10}). It follows from (\ref{eqInt.9})
that $\left\langle \,\cdot\,,\,\cdot\,\right\rangle _{\mathcal{H}_{s}}$ is
defined from the integral kernel $\left|  x-y\right|  ^{s-1}$. The
corresponding operator $A_{s}$ is a special case of the Knapp--Stein
intertwining operator, see \cite{KnSt80}. (See also \cite{Sal62} and
\cite{Rad98}.) This operator $A_{s}\left(  n\right)  $ is defined more
generally and also in $\mathbb{R}^{n}$. Then the integral kernel is $\left|
x-y\right|  ^{s-n}$, and $0<s<n$. If $\Delta$ is the positive Laplace operator
in $\mathbb{R}^{n}$, i.e., $\Delta=\sum_{j=1}^{n}\left(  \frac{1}{i}%
\frac{\partial\,}{\partial x_{j}}\right)  ^{2}$, then it is shown in
\cite[Lemma 2, p.~117]{Ste70} that $A_{s}=\Delta^{-\frac{s}{2}}$, and the
Fourier transform of $\left|  x\right|  ^{s-n}$ is
\[
\left(  \pi^{-\frac{s}{2}}\Gamma\left(  \frac{s}{2}\right)  \biggm/%
\Gamma\left(  \frac{n-s}{2}\right)  \right)  \cdot\left|  \xi\right|  ^{-s}.
\]
Hence up to a constant, the norm $\left\|  \,\cdot\,\right\|  _{s}$ of
(\ref{eqInt.9}) may be rewritten as%
\begin{equation}
\int_{\mathbb{R}}\left|  \xi\right|  ^{-s}\left|  \hat{f}\left(  \xi\right)
\right|  ^{2}\,d\xi,\label{eqRef.14}%
\end{equation}
and the inner product $\left\langle \,\cdot\,,\,\cdot\,\right\rangle _{s}$ as
\begin{equation}
\int_{\mathbb{R}}\left|  \xi\right|  ^{-s}\overline{\hat{f}_{1}\left(
\xi\right)  }\hat{f}_{2}\left(  \xi\right)  \,d\xi,\label{eqRef.15}%
\end{equation}
where
\begin{equation}
\hat{f}\left(  \xi\right)  =\int_{\mathbb{R}}e^{-i\xi x}f\left(  x\right)
\,dx\label{eqRef.16}%
\end{equation}
is the usual Fourier transform suitably extended to $\mathcal{H}_{s}$, using
Stein's singular integrals. Intuitively, $\mathcal{H}_{s}$ consists of
functions on $\mathbb{R}$ which arise as $\left(  \frac{d\,}{dx}\right)
^{s}f_{s}$ for some $f_{s}$ in $L^{2}\left(  \mathbb{R}\right)  $. This also
introduces a degree of ``non-locality'' into the theory, and the functions in
$\mathcal{H}_{s}$ cannot be viewed as locally integrable, although
$\mathcal{H}_{s}$ for each $s$, $0<s<1$, contains $C_{c}\left(  \mathbb{R}%
\right)  $ as a dense subspace. In fact, formula (\ref{eqRef.14}), for the
norm in $\mathcal{H}_{s}$, makes precise in which sense elements of
$\mathcal{H}_{s}$ are ``fractional'' derivatives of locally integrable
functions on $\mathbb{R}$, and that there are elements of $\mathcal{H}_{s}$
(and of $\mathcal{K}_{s}$) which are not locally integrable. On the other
hand, vectors in $\mathcal{H}_{s}$ are not too singular: for example the Dirac
function $\delta$ is not in $\mathcal{H}_{s}$. To see this, pick some
approximate identity $\varphi_{\varepsilon}\underset{\varepsilon\rightarrow
0}{\longrightarrow}\delta$, say $\varphi\in C_{c}\left(  -1,1\right)  $,
$\varphi>0$, $\int\varphi\left(  x\right)  \,dx=1$, and set $\varphi
_{\varepsilon}\left(  x\right)  =\frac{1}{\varepsilon}\varphi\left(  \frac
{x}{\varepsilon}\right)  $; then a calculation shows that%
\begin{equation}
\left\|  \varphi_{\varepsilon}\right\|  _{\mathcal{H}_{s}}^{2}=C_{s}%
\varepsilon^{s-1}\label{eqRef.17}%
\end{equation}
for some positive constant $C_{s}$. Hence $\delta$ is not in $\mathcal{H}_{s}%
$, and then of course also not in the subspace $\mathcal{K}_{s}$.

Nonetheless, if we pass to the new norm $\left\|  f\right\|  _{J}^{2}=\left\|
f\right\|  _{\mathcal{H}\left(  \mathcal{K}_{s}\right)  }^{2}=\left\langle
f,Jf\right\rangle _{s}$ of (\ref{eqRef.13}), then from (\ref{eqInt.10}) we get%
\begin{equation}
\left\|  \varphi_{\varepsilon}\right\|  _{J}^{2}=\mathcal{O}\left(
\varepsilon^{2}\right)  . \label{eqRef.18}%
\end{equation}
Hence the limit $\varphi_{\varepsilon}\rightarrow\delta$ defines a bounded
linear functional on $\mathcal{H}\left(  \mathcal{K}_{s}\right)  $, relative
to the norm $\left\|  \,\cdot\,\right\|  _{J}$ on that Hilbert space. From the
Riesz lemma, and the definition of $\mathcal{H}\left(  \mathcal{K}_{s}\right)
$, we conclude that $\delta$ is in $\mathcal{H}\left(  \mathcal{K}_{s}\right)
$. The same argument shows that the distributions $\delta^{\left(  n\right)
}:=\left(  \frac{d\,}{dx}\right)  ^{n}\delta$ given by
\begin{equation}
\delta^{\left(  n\right)  }\left(  \phi\right)  =\left(  -1\right)  ^{n}%
\frac{d^{n}\phi}{dx^{n}}\left(  0\right)  \label{eqRef.19}%
\end{equation}
for $\phi\in C_{c}^{\infty}\left(  -1,1\right)  $, are also in $\mathcal{H}%
\left(  \mathcal{K}_{s}\right)  $. In fact, the norm computes out as%
\begin{equation}
\left\|  \delta^{\left(  n\right)  }\right\|  _{J}^{2}=n!\left(  1-s\right)
\left(  2-s\right)  \cdots\left(  n-s\right)  \text{\qquad for }n=0,1,2,\dots.
\label{eqRef.20}%
\end{equation}

In the next lemma we provide the detailed proof of the fact that the iterated
derivatives $\left(  \frac{d\,}{dx}\right)  ^{n}\delta=:\delta^{\left(
n\right)  }$ of the Dirac delta function are all in the completion of
$C_{c}^{\infty}\left(  -1,1\right)  $ relative to the ``new'' norm of the
Hilbert space $\mathcal{H}\left(  \mathcal{K}_{s}\right)  $. But recall that
$\delta$, or its derivatives, are not in $\mathcal{K}_{s}$.

\begin{lemma}
\label{LemDiracsinHKs}For the Dirac mass and its derivatives, we have
$\delta^{\left(  n\right)  }\in\mathcal{H}\left(  \mathcal{K}_{s}\right)  $,
$n=0,1,2,\dots$. The restriction on $s$ is, as before, $0<s<1$.
\end{lemma}

\begin{proof}
First note that if $\phi\in C_{c}^{\infty}\left(  -1,1\right)  $, then
\begin{equation}
\int_{-1}^{1}\phi\left(  x\right)  \left(  1-xy\right)  ^{s-1}\,dx
\label{eqDiracsinHKs.1}%
\end{equation}
restricts to a $C^{\infty}$-function on $\left[  -1,1\right]  $. By this we
mean that there is a $C^{\infty}$-function $\varphi_{s}$ on $\mathbb{R}$ such
that%
\begin{equation}
\varphi_{s}\left(  y\right)  =\int_{-1}^{1}\phi\left(  x\right)  \left(
1-xy\right)  ^{s-1}\,dx \label{eqDiracsinHKs.2}%
\end{equation}
holds for all $y$ in $\left[  -1,1\right]  $. Hence, if $F$ is a distribution
with compact support in $\left[  -1,1\right]  $, then%
\begin{equation}
\left\langle \varphi_{s},F\right\rangle =F\left(  \varphi_{s}\right)
\label{eqDiracsinHKs.3}%
\end{equation}
is well-defined. The same argument shows that $\left\langle \left(
1-\,\cdot\;y\right)  ^{s-1},F\right\rangle $ is well-defined, and that%
\[
y\longmapsto\left\langle \left(  1-\,\cdot\;y\right)  ^{s-1},F\right\rangle
\]
is also $C^{\infty}$ up to the endpoints in the closed interval $I=\left[
-1,1\right]  $. Hence, the distribution $F$ may be applied again, and we get
the expression%
\begin{equation}
\left\|  F\right\|  _{\mathcal{H}\left(  \mathcal{K}_{s}\right)  }^{2}%
:=\int_{I}\int_{I}\overline{F\left(  x\right)  }\left(  1-xy\right)
^{s-1}F\left(  y\right)  \,dx\,dy. \label{eqDiracsinHKs.4}%
\end{equation}
Moreover, if $\phi\in C_{c}^{\infty}\left(  -1,1\right)  $, then
\[
\left\langle W\phi,F\right\rangle _{\mathcal{H}\left(  \mathcal{K}_{s}\right)
}=\int_{I}\int_{I}\overline{\phi\left(  x\right)  }\left(  1-xy\right)
^{s-1}F\left(  y\right)  \,dx\,dy
\]
is well-defined in the distribution sense, and%
\[
\left|  \left\langle W\phi,F\right\rangle _{\mathcal{H}\left(  \mathcal{K}%
_{s}\right)  }\right|  \leq\left\|  W\phi\right\|  _{\mathcal{H}\left(
\mathcal{K}_{s}\right)  }\left\|  F\right\|  _{\mathcal{H}\left(
\mathcal{K}_{s}\right)  },
\]
where $\left\|  F\right\|  _{\mathcal{H}\left(  \mathcal{K}_{s}\right)  }$ is
the expression (\ref{eqDiracsinHKs.4}). Hence for each $n=0,1,2,\dots$, we
must show the following implication:%
\begin{equation}
\left\langle W\phi,F\right\rangle _{\mathcal{H}\left(  \mathcal{K}_{s}\right)
}=0\text{ for all }\phi\in C_{c}^{\infty}\left(  -1,1\right)  \Longrightarrow
\left\langle \delta^{\left(  n\right)  },F\right\rangle _{\mathcal{H}\left(
\mathcal{K}_{s}\right)  }=0. \label{eqDiracsinHKs.5}%
\end{equation}
The interpretation of the brackets $\left\langle \,\cdot\,,\,\cdot
\,\right\rangle _{\mathcal{H}\left(  \mathcal{K}_{s}\right)  }$ is in the
sense of distributions as noted. In particular,
\begin{equation}
\left\langle \delta^{\left(  n\right)  },F\right\rangle _{\mathcal{H}\left(
\mathcal{K}_{s}\right)  }=\left(  s-1\right)  \cdots\left(  s-n\right)
\int_{I}y^{n}F\left(  y\right)  \,dy, \label{eqDiracsinHKs.6}%
\end{equation}
where $\int_{I}y^{n}F\left(  y\right)  \,dy$ is really the compactly supported
distribution $F$ evaluated at the monomial $y^{n}$. Recall, it is assumed that
the distribution $F$ is supported in $I$. Now pick $\phi\in C_{c}^{\infty
}\left(  -1,1\right)  $ such that $\phi>0$, and $\int_{I}\phi\left(  x\right)
\,dx=1$, and let $\phi_{\varepsilon}\left(  x\right)  =\frac{1}{\varepsilon
}\phi\left(  \frac{x}{\varepsilon}\right)  $, for $0<\varepsilon<1$. We prove
next that
\begin{equation}
\lim_{\varepsilon\rightarrow0}\left\langle W\phi_{\varepsilon}^{\left(
n\right)  },F\right\rangle _{\mathcal{H}\left(  \mathcal{K}_{s}\right)
}=\left\langle \delta^{\left(  n\right)  },F\right\rangle _{\mathcal{H}\left(
\mathcal{K}_{s}\right)  }, \label{eqDiracsinHKs.7}%
\end{equation}
where both sides are understood in the sense of distributions. But we also
have $\left\langle W\phi_{\varepsilon}^{\left(  n\right)  },F\right\rangle =0$
for all $\varepsilon>0$, by the assumption in (\ref{eqDiracsinHKs.5}). To
complete the proof we will then only need to check that%
\begin{equation}
\sup_{0<\varepsilon<1}\left\|  W\phi_{\varepsilon}^{\left(  n\right)
}\right\|  _{\mathcal{H}\left(  \mathcal{K}_{s}\right)  }<\infty.
\label{eqDiracsinHKs.8}%
\end{equation}
Explicitly,%
\begin{equation}
\left\|  W\phi_{\varepsilon}^{\left(  n\right)  }\right\|  _{\mathcal{H}%
\left(  \mathcal{K}_{s}\right)  }^{2}=\int_{I}\int_{I}\phi_{\varepsilon
}^{\left(  n\right)  }\left(  x\right)  \left(  1-xy\right)  ^{s-1}%
\phi_{\varepsilon}^{\left(  n\right)  }\left(  y\right)  \,dx\,dy,
\label{eqDiracsinHKs.9}%
\end{equation}
and this last expression can be estimated directly: If $n\in\left\{
0,1,2,\dots\right\}  $, there is a constant $C_{n}$ ($<\infty$) such that the
$\int_{I}\int_{I}\cdots\,dx\,dy$ term in (\ref{eqDiracsinHKs.9}) is estimated
by $C_{n}$. In particular, we have the desired estimate (\ref{eqDiracsinHKs.8}%
). The left-hand side of (\ref{eqDiracsinHKs.7}) may therefore be estimated by
$\left\|  F\right\|  _{\mathcal{H}\left(  \mathcal{K}_{s}\right)  }\cdot
C_{n}$. Since $\left\langle W\phi_{\varepsilon}^{\left(  n\right)
},F\right\rangle _{\mathcal{H}\left(  \mathcal{K}_{s}\right)  }=0$ for all $n$
and all $\varepsilon$, by assumption, see (\ref{eqDiracsinHKs.5}), we will
then have $\left\langle \delta^{\left(  n\right)  },F\right\rangle
_{\mathcal{H}\left(  \mathcal{K}_{s}\right)  }=0$, which is the claim.

It remains to check that the limit (as $\varepsilon\rightarrow0$) in
(\ref{eqDiracsinHKs.7}) is as stated. The argument is much as the previous
one, so we will merely sketch the details for the case of $n=0$: Since $F$ is
an distribution with support in $I=\left[  -1,1\right]  $, we need to check
that
\begin{equation}
\lim_{\varepsilon\rightarrow0}\frac{1}{\varepsilon}\int_{I}\phi\left(
\frac{x}{\varepsilon}\right)  \left(  1-xy\right)  ^{s-1}\,dx=1
\label{eqDiracsinHKs.10}%
\end{equation}
and%
\begin{equation}
\lim_{\varepsilon\rightarrow0}\frac{1}{\varepsilon}\left(  \frac{d\,}%
{dy}\right)  ^{m}\int_{I}\phi\left(  \frac{x}{\varepsilon}\right)  \left(
1-xy\right)  ^{s-1}\,dx=0 \label{eqDiracsinHKs.11}%
\end{equation}
for all $m\in\mathbb{N}$; and both of these limits can be verified by
calculus. Indeed the left-hand side in (\ref{eqDiracsinHKs.10}) is of the
order%
\[
L_{\varepsilon}^{\left(  s\right)  }\left(  y\right)  :=%
%TCIMACRO{\TeXButton{cases}{\begin{cases}
%\displaystyle\frac{\left( 1+\varepsilon y\right) ^{s}-\left( 1-\varepsilon
%y\right) ^{s}}{2\varepsilon sy}&\text{if }y\neq0,\\1&\text{if }y=0,
%\end{cases}
%}}%
%BeginExpansion
\begin{cases}
\displaystyle\frac{\left( 1+\varepsilon y\right) ^{s}-\left( 1-\varepsilon
y\right) ^{s}}{2\varepsilon sy}&\text{if }y\neq0,\\1&\text{if }y=0,
\end{cases}
%EndExpansion
\]
which is differentiable in $y$, for every $\varepsilon\in\mathbb{R}_{+}$. The
corresponding expression in (\ref{eqDiracsinHKs.11}) is $\mathcal{O}\left(
\varepsilon^{m}\right)  $, $m=1,2,\dots$. Since the distribution is of compact
support (in $I$) we also have, for some $m\in\mathbb{N}$, the estimate
\[
\left|  F\left(  \psi\right)  \right|  \leq\volapuk\cdot\max_{0\leq k\leq
m}\max_{x\in I}\left|  \psi^{\left(  k\right)  }\left(  x\right)  \right|
\]
for all $\psi\in C^{\infty}\left(  \mathbb{R}\right)  $.

Applying this to the functions $\psi$ ($=L_{\varepsilon}$) in the left-hand
side of (\ref{eqDiracsinHKs.10}), we finally arrive at the desired conclusion
(\ref{eqDiracsinHKs.7}). This completes the proof of the lemma.
\end{proof}

Hence if $f\in\mathcal{K}_{s}$, $Wf\in\mathcal{H}\left(  \mathcal{K}%
_{s}\right)  $, we get the inner product $\left\langle \delta^{\left(
n\right)  },Wf\right\rangle _{J}$ is well-defined. A calculation yields%
\begin{equation}
\left\langle \delta^{\left(  n\right)  },Wf\right\rangle _{J}=\left(
s-1\right)  \left(  s-2\right)  \cdots\left(  s-n\right)  \int_{-1}^{1}%
x^{n}f\left(  x\right)  \,dx. \label{eqRef.21}%
\end{equation}
However, if $f$ is not locally integrable, then the right-hand side in
(\ref{eqRef.21}) must be understood as a singular integral, see, e.g.,
\cite[Chapters V.1--2]{Ste70}.

Recall that $\mathcal{K}_{s}$ is obtained as the completion of $C_{c}\left(
-1,1\right)  $ relative to the norm $\left\|  \,\cdot\,\right\|  _{s}$ of
(\ref{eqInt.9}). If $f$ is in $C_{c}\left(  -1,1\right)  $, then the Fourier
transform%
\begin{equation}
\hat{f}\left(  \xi\right)  =\int_{-1}^{1}f\left(  x\right)  e^{-ix\xi}\,dx
\label{eqRef.22}%
\end{equation}
of (\ref{eqRef.16}) clearly has an entire analytic extension, i.e., it extends
to complex values of $\xi$ as an entire analytic function with exponential
growth factor $e^{\left|  \operatorname{Im}\xi\right|  }$, $\xi\in\mathbb{C}$.
We wish to show that this also holds for $f\in\mathcal{N}\subset
\mathcal{K}_{s}$. Note if $f\in\mathcal{N}$, it has finite $\left\|
\,\cdot\,\right\|  _{s}$-norm, and%
\begin{equation}
\int_{-1}^{1}\int_{-1}^{1}\overline{f\left(  x\right)  }\left(  1-xy\right)
^{s-1}f\left(  y\right)  \,dx\,dy=0, \label{eqRef.23}%
\end{equation}
or rather $\left\|  f\right\|  _{J}=0$. Since $f$ can be rather singular, the
claim requires a proof. We have $Wf=0$, and the Dirac measures $\delta_{x}$,
for $x\in\mathbb{R}$, $\left|  x\right|  <1$, are in $\mathcal{H}\left(
\mathcal{K}_{s}\right)  $. Hence $\left\langle \delta_{x},Wf\right\rangle
_{J}=0$. But a calculation yields, for $x\in\left(  -1,1\right)  =:I$,%
\begin{equation}
\left\langle \delta_{x},Wf\right\rangle _{J}=\int_{-1}^{1}\left(  1-xy\right)
^{s-1}f\left(  y\right)  \,dy. \label{eqRef.24}%
\end{equation}
Let $x\in I\setminus\left\{  0\right\}  $, and multiply by $\left|  x\right|
^{1-s}$, to get%
\[
\int_{-1}^{1}\left|  \frac{1}{x}-y\right|  ^{s-1}f\left(  y\right)  \,dy=0,
\]
and so $\left(  A_{s}f\right)  \left(  \frac{1}{x}\right)  =0$. We conclude
that $Af$ is supported in the interval if $f$ is in $\mathcal{N}$. This
localizes the computation of
\begin{equation}
\left\|  f\right\|  _{s}^{2}=\int_{\mathbb{R}}\overline{f\left(  x\right)
}A_{s}f\left(  x\right)  \,dx, \label{eqRef.25}%
\end{equation}
but still interpreted as a singular integral.

Since $\left\|  f\right\|  _{s}<\infty$, and $f\in\mathcal{K}_{s}$, there is a
sequence $\varphi_{n}\in C_{c}^{\infty}\left(  -1,1\right)  $ such that
$\lim_{n\rightarrow\infty}\left\|  f-\varphi_{n}\right\|  _{s}=0$. Then of
course also%
\begin{equation}
\lim_{n\rightarrow\infty}\left\|  \varphi_{n}\right\|  _{s}=\left\|
f\right\|  _{s}<\infty. \label{eqRef.26}%
\end{equation}
But%
\begin{equation}
\left\|  \varphi_{n}\right\|  _{s}^{2}=C_{s}\int_{\mathbb{R}}\left|
\xi\right|  ^{-s}\left|  \hat{\varphi}_{n}\left(  \xi\right)  \right|
^{2}\,d\xi\label{eqRef.27}%
\end{equation}
by (\ref{eqRef.14}). It follows that there is a subsequence $\varphi_{n_{i}}$
such that $\hat{\varphi}_{n_{i}}\left(  \,\cdot\,\right)  $ converges
pointwise almost everywhere on $\mathbb{R}$. We wish to use Montel's theorem
\cite[v.~II, Theorem 15.3.1]{Hil62} to conclude that the Fourier transform
$\hat{f}$ of $f$ also has an entire analytic extension. To do this we need
only check that $\hat{\varphi}_{n_{i}}\left(  \zeta\right)  $, $\zeta
\in\mathbb{C}$, is an equicontinuous family. Now pick $\zeta_{1},\zeta_{2}%
\in\mathbb{C}$, and consider%
\[
\hat{\varphi}_{n_{i}}\left(  \zeta_{1}\right)  -\hat{\varphi}_{n_{i}}\left(
\zeta_{2}\right)  =\int_{-1}^{1}\varphi_{n_{i}}\left(  x\right)  \left\{
e^{-ix\zeta_{1}}-e^{-ix\zeta_{2}}\right\}  \,dx.
\]
Let $E\left(  x\right)  :=e^{-ix\zeta_{1}}-e^{-ix\zeta_{2}}$, and pick
$\psi\in C_{c}^{\infty}\left(  \mathbb{R}\right)  $ such that $\psi\equiv1$ on
$\bar{I}=\left[  -1,1\right]  $. Continuing the calculation, we get%
\begin{align*}
\int_{-1}^{1}\varphi_{n_{i}}\left(  x\right)  E\left(  x\right)  \,dx  &
=\int_{\mathbb{R}}\varphi_{n_{i}}\left(  x\right)  \psi\left(  x\right)
E\left(  x\right)  \,dx\\
&  =\int_{\mathbb{R}}\left(  \Delta^{-\frac{s}{2}}\varphi_{n_{i}}\left(
x\right)  \right)  \left(  \Delta^{\frac{s}{2}}\psi E\left(  x\right)
\right)  \,dx
\end{align*}
and%
\begin{align*}
\left|  \int_{-1}^{1}\varphi_{n_{i}}\left(  x\right)  E\left(  x\right)
\,dx\right|   &  \leq\left\|  \varphi_{n_{i}}\right\|  _{s}\cdot\left\|
\Delta^{\frac{s}{2}}\psi E\right\|  _{L^{2}\left(  \mathbb{R}\right)  }\\
&  \leq\left\|  \varphi_{n_{i}}\right\|  _{s}\cdot\left\{  \int_{\mathbb{R}%
}\left|  \frac{d\,}{dx}\left(  \psi E\right)  \left(  x\right)  \right|
^{2}\,dx\right\}  ^{\frac{1}{2}}.
\end{align*}
But we have from (\ref{eqRef.26}) that $\sup_{i}\left\|  \varphi_{n_{i}%
}\right\|  _{s}<\infty$, and the second term is independent of $n_{i}$, and it
can be estimated in terms of $\left|  \zeta_{1}-\zeta_{2}\right|  $ by
calculus. This shows that the entire functions $\left\{  \hat{\varphi}_{n_{i}%
}\left(  \zeta\right)  \right\}  $ do form an equicontinuous family. Since
$\hat{\varphi}_{n_{i}}\left(  \xi\right)  $ is convergent $\mathrm{a.e.}%
\;\xi\in\mathbb{R}$ as noted, we conclude that the entire functions
$\hat{\varphi}_{n_{i}}\left(  \zeta\right)  $ converge uniformly for $\zeta$
in compact subsets of $\mathbb{C}$, and that the limit function is also entire
analytic. But by the argument above, this limit is an extension of $\hat
{f}\left(  \xi\right)  $, for $\xi\in\mathbb{R}$. From (\ref{eqRef.21}), we
have%
\[
\left\langle \delta^{\left(  n\right)  },Wf\right\rangle _{J}=\left(
s-1\right)  \left(  s-2\right)  \cdots\left(  s-n\right)  i^{n}\left(
\frac{d\,}{d\zeta}\right)  ^{n}\hat{f}\left(  \zeta\right)  |_{\zeta=0}.
\]
Since $f\in\mathcal{N}$, $Wf=0$, and the left-hand side vanishes for all
$n=0,1,2,\dots$. Hence all the derivatives $\left(  \frac{d\,}{d\zeta}\right)
^{n}\hat{f}\left(  \zeta\right)  $ vanish at $\zeta=0$. Since $\hat{f}$ is
analytic, it must vanish identically. Finally use (\ref{eqRef.14}) to conclude
that $f=0$ as an element of $\mathcal{K}_{s}$. This completes the proof of
(\ref{ThmRef.1(3)}), and therefore the proof of the theorem.
\end{proof}

In Section \ref{Han}, we will consider more systematically the structure of
systems $\left(  \mathcal{H}_{0},\mathcal{K},J,U\right)  $ for which
$W\colon\mathcal{K}\rightarrow\mathcal{H}\left(  \mathcal{K}\right)  $ is
$1$--$1$. The present construction (i.e., Theorem \ref{ThmRef.1}%
(\ref{ThmRef.1(3)})) has the initial operator $U$ unitary in $\mathcal{H}_{0}%
$, and in fact part of a unitary one-parameter group. If the unitarity
restriction on $U$ is relaxed, then there is a richer variety of examples with
$\ker\left(  W\right)  =\left\{  0\right\}  $. For example, let $A$ denote the
\emph{unilateral shift} in $H^{2}=H^{2}\left(  \mathbb{T}\right)  $, and set
\[
U=%
\begin{pmatrix}
A & 0\\
0 & A^{\ast}%
\end{pmatrix}
,\qquad J=%
\begin{pmatrix}
0 & I\\
I & 0
\end{pmatrix}
\]
on $\mathcal{H}_{0}=H^{2}\oplus H^{2}$. Then we show in Section \ref{Han} that
the subspaces $\mathcal{K}$ described axiomatically in Theorem \ref{ThmRef.1}
above, and which are further assumed maximal, are in $1$--$1$ correspondence
with finite positive Borel measures on $\left[  -1,1\right]  $, such that
$n\mapsto\int x^{n}\,d\mu\left(  x\right)  $ is in $\ell^{2}$. For those
examples, the condition $\ker\left(  W_{\mu}\right)  =\left\{  0\right\}  $
holds if and only if $\operatorname*{supp}\left(  \mu\right)  $ has
accumulation points in $\left(  -1,1\right)  $. It holds, for example, if
$\mu$ is the restriction to $\left[  -1,1\right]  $ of Lebesgue measure.

\section{\label{Rep}Reproducing kernels}

In the proof of part (\ref{ThmRef.1(3)}) of Theorem \ref{ThmRef.1}, we used
the reflection $J$ to arrive at a new Hilbert space $\mathcal{H}\left(
\mathcal{K}_{s}\right)  $. Recall that $\mathcal{K}_{s}$ is the closure of
$C_{c}\left(  -1,1\right)  $ in the norm $\left\|  \,\cdot\,\right\|  _{s}$
defined as in (\ref{eqInt.9}) from the Knapp--Stein operator $A_{s}$. But in
part (\ref{ThmRef.1(2)}) of Theorem \ref{ThmRef.1}, we showed that the system
$\left(  \mathcal{H}\left(  \mathcal{K}_{s}\right)  ,W,S\right)  $ is
determined uniquely from $\left(  \mathcal{K}_{s},J\right)  $ up to unitary
equivalence. In proving part (\ref{ThmRef.1(3)}), we selected a particular
version of $\mathcal{H}\left(  \mathcal{K}_{s}\right)  $ which turned out to
contain distributions, specifically, we showed that $\left\{  \delta^{\left(
n\right)  }=\left(  \frac{d\,}{dx}\right)  ^{n}\delta\mathrel{;}%
n=0,1,\dots\right\}  $ forms an orthogonal basis in $\mathcal{H}\left(
\mathcal{K}_{s}\right)  $. Our interpretation of this is that we make the
Taylor expansion around $x=0$ into an orthogonal expansion relative to the
inner product in $\mathcal{H}\left(  \mathcal{K}_{s}\right)  $. But there is
an alternative construction of $\mathcal{H}\left(  \mathcal{K}_{s}\right)  $
consisting of analytic functions in
\begin{equation}
D:=\left\{  z\in\mathbb{C}\mathrel{;}\left|  z\right|  <1\right\}  .
\label{eqRep.1}%
\end{equation}
This is a Hilbert space $\mathcal{H}_{\mathrm{rep}}\left(  s\right)  $
constructed as a reproducing kernel Hilbert space from the kernel%
\begin{equation}
Q_{s}\left(  z,w\right)  =\left(  1-z\bar{w}\right)  ^{s-1},\qquad\left(
z,w\right)  \in D\times D. \label{eqRep.2}%
\end{equation}
It is known that there is a unique Hilbert space $\mathcal{H}_{\mathrm{rep}%
}\left(  s\right)  $ consisting of analytic functions on $D$ such that%
\begin{equation}
f\left(  w\right)  =\left\langle Q_{s}\left(  \,\cdot\,,w\right)
,f\right\rangle _{\mathcal{H}_{\mathrm{rep}}\left(  s\right)  },
\label{eqRep.3}%
\end{equation}
where $\left\langle \,\cdot\,,\,\cdot\,\right\rangle _{\mathcal{H}%
_{\mathrm{rep}}\left(  s\right)  }$ is the inner product of this Hilbert
space. It has the monomials $\left\{  z^{n}\mathrel{;}n=0,1,2,\dots\right\}  $
as an orthogonal basis, and we refer to \cite{ShSh62} and \cite{Aro50} for
more details on these Hilbert spaces. It will be convenient for us to denote
the kernel functions in $\mathcal{H}_{\mathrm{rep}}\left(  s\right)  $,
\begin{equation}
q_{w}\left(  z\right)  :=\left(  1-\bar{w}z\right)  ^{s-1}. \label{eqRep.4}%
\end{equation}
An application of (\ref{eqRep.3}) then yields%
\begin{equation}
Q_{s}\left(  w_{1},w_{2}\right)  =\left\langle q_{w_{1}},q_{w_{2}%
}\right\rangle _{\mathcal{H}_{\mathrm{rep}}\left(  s\right)  }.
\label{eqRep.5}%
\end{equation}

\begin{corollary}
\label{CorRep.1}The two Hilbert spaces $\mathcal{H}\left(  \mathcal{K}%
_{s}\right)  $ and $\mathcal{H}_{\mathrm{rep}}\left(  s\right)  $, $0<s<1$,
are naturally isomorphic with a unitary isomorphism%
\begin{equation}
T\colon\mathcal{H}\left(  \mathcal{K}_{s}\right)  \longrightarrow
\mathcal{H}_{\mathrm{rep}}\left(  s\right)  \label{eqRep.6}%
\end{equation}
which intertwines the respective selfadjoint scaling operators%
\begin{equation}
\left(  S_{a}f\right)  \left(  x\right)  =a^{s+1}f\left(  a^{2}x\right)
\label{eqRep.7}%
\end{equation}
and%
\begin{equation}
\left(  S_{a}^{\mathbb{C}}F\right)  \left(  z\right)  =a^{s-1}F\left(
a^{-2}z\right)  , \label{eqRep.8}%
\end{equation}
for $f\in\mathcal{H}\left(  \mathcal{K}_{s}\right)  $, $x\in\mathbb{R}$,
$F\in\mathcal{H}_{\mathrm{rep}}\left(  s\right)  $, $z\in D$, $a>1$.
Specifically, we have
\begin{equation}
TS_{a}^{{}}=S_{a}^{\mathbb{C}}T. \label{eqRep.9}%
\end{equation}
\end{corollary}

\begin{proof}
While it is possible to give a direct proof along the lines of the last two
pages in section 9 of \cite{JoOl99}, we will derive the result here as a
direct corollary to Theorem \ref{ThmRef.1}(\ref{ThmRef.1(2)}), i.e., the
uniqueness up to unitary equivalence. Given $a>1$, we already established the
system $\left(  \mathcal{H}\left(  \mathcal{K}_{s}\right)  ,W,S_{a}\right)  $
in part (\ref{ThmRef.1(3)}) of Theorem \ref{ThmRef.1}. We wish to show that
there is a second system
\begin{equation}
\left(  \mathcal{H}_{\mathrm{rep}}^{{}}\left(  s\right)  ,W_{{}}^{\mathbb{C}%
},S_{a}^{\mathbb{C}}\right)  ,\qquad W_{{}}^{\mathbb{C}}=W_{s}^{\mathbb{C}},
\label{eqRep.10}%
\end{equation}
which also satisfies axioms (\ref{eqThmRef.1(3)})--(\ref{eqThmRef.1(4)}) in
part (\ref{ThmRef.1(2)}). The $s$-dependence of $W=W_{s}$ will be suppressed
in the proof for simplicity. For $S_{a}^{\mathbb{C}}$ we take the
transformation defined in (\ref{eqRep.8}) above, and we get $W^{\mathbb{C}%
}\colon\mathcal{K}_{s}\rightarrow\mathcal{H}_{\mathrm{rep}}\left(  s\right)  $
by the following formula:%
\begin{equation}
\left(  W^{\mathbb{C}}k\right)  \left(  z\right)  =\int_{-1}^{1}k\left(
x\right)  \left(  1-xz\right)  ^{s-1}\,dx \label{eqRep.11}%
\end{equation}
for $k\in\mathcal{K}_{s}$, and $z\in D$. To see that $S_{a}^{\mathbb{C}}$ in
(\ref{eqRep.8}) is selfadjoint in $\mathcal{H}_{\mathrm{rep}}\left(  s\right)
$, we compute the inner products as follows:%
\begin{align*}
\left\langle S_{a}^{\mathbb{C}}q_{w_{1}}^{{}},q_{w_{2}}^{{}}\right\rangle
_{\mathrm{rep}}  &  =a^{s-1}\left\langle q_{w_{1}}\left(  a^{-2}%
\,\cdot\,\right)  ,q_{w_{2}}\right\rangle _{\mathrm{rep}}\\
&  =a^{s-1}\left\langle q_{a^{-2}w_{1}}\left(  \,\cdot\,\right)  ,q_{w_{2}%
}\right\rangle _{\mathrm{rep}}\\
&  =a^{s-1}Q_{s}\left(  a^{-2}w_{1},w_{2}\right) \\
&  =a^{s-1}\left(  1-a^{-2}w_{1}\bar{w}_{2}\right)  ^{s-1}\\
&  =a^{s-1}Q_{s}\left(  w_{1},a^{-2}w_{2}\right) \\
&  =\left\langle q_{w_{1}}^{{}},S_{a}^{\mathbb{C}}q_{w_{2}}^{{}}\right\rangle
_{\mathrm{rep}}\text{\qquad for all }w_{1},w_{2}\in D.
\end{align*}
Since the kernel functions $\left\{  q_{w}^{\left(  s\right)  }\mathrel{;}w\in
D\right\}  $ are dense in $\mathcal{H}_{\mathrm{rep}}\left(  s\right)  $ by
construction, we conclude that $S_{a}^{\mathbb{C}}$ is indeed selfadjoint in
$\mathcal{H}_{\mathrm{rep}}\left(  s\right)  $ when $a>1$ and $0<s<1$.

We now show that $W^{\mathbb{C}}\colon\mathcal{K}_{s}\rightarrow
\mathcal{H}_{\mathrm{rep}}\left(  s\right)  $ in (\ref{eqRep.11}) is
contractive. For $k\in\mathcal{K}_{s}$, we have
\begin{align*}
\left\|  W^{\mathbb{C}}k\right\|  _{\mathrm{rep}}^{2}  &  =\int_{-1}^{1}%
\int_{-1}^{1}\overline{k\left(  x\right)  }\left\langle q_{x},q_{y}%
\right\rangle _{\mathrm{rep}}k\left(  y\right)  \,dx\,dy\\
&  =\int_{-1}^{1}\int_{-1}^{1}\overline{k\left(  x\right)  }\left(
1-xy\right)  ^{s-1}k\left(  y\right)  \,dx\,dy\\
&  =\int_{\mathbb{R}}\overline{k\left(  x\right)  }A_{s}Jk\left(  x\right)
\,dx\\
&  =\left\langle k,Jk\right\rangle _{\mathcal{H}_{s}}\leq\left\|  k\right\|
_{s}^{2},
\end{align*}
which shows that $W_{s}^{\mathbb{C}}$ is contractive as claimed. But we also
proved that%
\[
\left\langle W^{\mathbb{C}}k_{1},W^{\mathbb{C}}k_{2}\right\rangle
_{\mathrm{rep}}=\left\langle k_{1},Jk_{2}\right\rangle _{\mathcal{H}_{s}}%
\]
for all $k_{1},k_{2}\in\mathcal{K}_{s}\subset\mathcal{H}_{s}$. Hence%
\begin{equation}
\left(  W^{\mathbb{C}}\right)  ^{\ast}W^{\mathbb{C}}=P_{s}JP_{s},
\label{eqRep.12}%
\end{equation}
where $P_{s}$ denotes the projection of $\mathcal{H}_{s}$ onto $\mathcal{K}%
_{s}$. Hence axiom (\ref{eqThmRef.1(4)}) in the statement of Theorem
\ref{ThmRef.1}(\ref{ThmRef.1(2)}) is also satisfied. We leave the verification
of%
\begin{equation}
S_{a}^{\mathbb{C}}W_{{}}^{\mathbb{C}}=W_{{}}^{\mathbb{C}}UP_{s}^{{}}
\label{eqRep.13}%
\end{equation}
from (\ref{ThmRef.1(2)})(\ref{eqThmRef.1(3)}) to the reader. The conclusion of
Corollary \ref{CorRep.1} is now immediate from Theorem \ref{ThmRef.1}%
(\ref{ThmRef.1(2)}).
\end{proof}

Let $T\colon\mathcal{H}\left(  \mathcal{K}_{s}\right)  \rightarrow
\mathcal{H}_{\mathrm{rep}}\left(  s\right)  $ be the unitary isomorphism from
(\ref{eqRep.9}) in the statement of Corollary \ref{CorRep.1}. We saw in
Theorem \ref{ThmRef.1}(\ref{ThmRef.1(2)}) that%
\[
TW_{s}^{{}}=W_{s}^{\mathbb{C}}T.
\]
Recall that $\delta^{\left(  n\right)  }=\left(  \frac{d\,}{dx}\right)
^{n}\delta$ is in $\mathcal{H}\left(  \mathcal{K}_{s}\right)  $, and we
conclude that
\[
T\left(  \delta^{\left(  n\right)  }\right)  \left(  z\right)  =\left(
s-1\right)  \left(  s-2\right)  \cdots\left(  s-n\right)  z^{n}.
\]
Since $T$ is isometric, and
\[
\left\|  \delta^{\left(  n\right)  }\right\|  _{\mathcal{H}\left(
\mathcal{K}_{s}\right)  }^{2}=\left(  1-s\right)  \cdots\left(  n-s\right)
n!\,,
\]
we conclude that%
\[
\left\|  z^{n}\right\|  _{\mathcal{H}_{\mathrm{rep}}\left(  s\right)  }%
^{2}=\frac{n!}{\left(  1-s\right)  \left(  2-s\right)  \cdots\left(
n-s\right)  }.
\]
We have proved the following

\begin{corollary}
\label{CorRep.2}Elements of $\mathcal{H}_{\mathrm{rep}}\left(  s\right)  $ may
be characterized by the orthogonal expansion%
\begin{align*}
f\left(  z\right)   &  =\sum_{n=0}^{\infty}c_{n}z^{n},\\
\left\|  f\right\|  _{\mathcal{H}_{\mathrm{rep}}\left(  s\right)  }^{2}  &
=\sum_{n=0}^{\infty}\left|  c_{n}\right|  ^{2}\frac{n!}{\left(  1-s\right)
\left(  2-s\right)  \cdots\left(  n-s\right)  }.
\end{align*}
\end{corollary}

\section{\label{Har}The Hardy space $H^{2}\left(  \mathbb{T}\right)  $}

In this section, we return to the space $L^{2}\left(  \mathbb{T}\right)  $ and
its subspace $H^{2}\left(  \mathbb{T}\right)  $ introduced in Section
\ref{Int}. Relative to the reflection $Jf\left(  z\right)  =f\left(  \bar
{z}\right)  $, $f\in L^{2}\left(  \mathbb{T}\right)  $, we describe a family
of positive subspaces defined from $H^{2}\left(  \mathbb{T}\right)  $. The
individual subspaces $\mathcal{K}\left(  b\right)  $ are positive relative to
$J$ and indexed by some function, $b$, say, in $H^{\infty}\left(
\mathbb{T}\right)  $. However, unless $b\equiv1$, the subspace $\mathcal{K}%
\left(  b\right)  $ is not shift invariant.

We first return to the axiomatic setup from Section \ref{Int}, and we derive a
formula for the contractive operator%
\begin{equation}
W\colon\mathcal{K}\longrightarrow\mathcal{H}\left(  \mathcal{K}\right)
\label{eqHar.1}%
\end{equation}
constructed from a given positive subspace $\mathcal{K}\subset\mathcal{H}_{0}%
$. Let $\mathcal{H}_{0}$ be a Hilbert space, and let $J$ be a period-$2$
unitary operator in $\mathcal{H}_{0}$. Let $\mathcal{H}_{\pm}$ be the
$J$-eigenspaces corresponding to eigenvalues $\pm1$, and let $P_{\pm}$ be the
respective projections onto $\mathcal{H}_{\pm}$, specifically%
\begin{equation}
P_{\pm}=\frac{1}{2}\left(  I\pm J\right)  . \label{eqHar.1prime}%
\end{equation}
We say that a closed subspace $\mathcal{K}\subset\mathcal{H}_{0}$ is
\emph{positive} if%
\begin{equation}
\left\langle k,Jk\right\rangle \geq0\text{\qquad for all }k\in\mathcal{K}.
\label{eqHar.2}%
\end{equation}
In Section \ref{Int}, we proved the following:

\begin{lemma}
\label{LemHar.1}\renewcommand{\theenumi}{\alph{enumi}}\renewcommand
{\theenumii}{\roman{enumii}}\renewcommand{\labelenumii}{(\theenumii)}

\begin{enumerate}
\item \label{LemHar.1(1)}There is a $1$--$1$ correspondence between the
following data \textup{(\ref{LemHar.1(1)(1)})} and
\textup{(\ref{LemHar.1(1)(2)}):}

\begin{enumerate}
\item \label{LemHar.1(1)(1)}closed positive subspaces $\mathcal{K}$,
\end{enumerate}

and

\begin{enumerate}
\stepcounter{enumii}

\item \label{LemHar.1(1)(2)}closed subspaces $\mathcal{K}_{+}\subset
\mathcal{H}_{+}$, and contractive linear operators%
\begin{equation}
\Lambda\colon\mathcal{K}_{+}\longrightarrow\mathcal{H}_{-}. \label{eqHar.3}%
\end{equation}
\end{enumerate}

\item \label{LemHar.1(2)}Given \textup{(\ref{LemHar.1(1)(1)}),} set%
\begin{equation}
\mathcal{K}_{+}:=P_{+}\mathcal{K}, \label{eqHar.4}%
\end{equation}
and%
\begin{equation}
\Lambda\left(  P_{+}k\right)  :=P_{-}k\text{\qquad for }k\in\mathcal{K}.
\label{eqHar.5}%
\end{equation}

\item \label{LemHar.1(3)}Given \textup{(\ref{LemHar.1(1)(2)}),} set
$\mathcal{K}:=G\left(  \Lambda\right)  =$ the graph of the contraction
$\Lambda$ in \textup{(\ref{eqHar.3}),} i.e.,%
\begin{equation}
\mathcal{K=}\left\{  k_{+}\oplus\Lambda k_{+}\mathrel{;}k_{+}\in
\mathcal{K}_{+}\right\}  . \label{eqHar.6}%
\end{equation}
\end{enumerate}
\end{lemma}

\begin{proof}
While the details are essentially in Section \ref{Int}, we sketch
(\ref{LemHar.1(1)(1)}) $\leftrightarrow$ (\ref{LemHar.1(1)(2)}).
(\ref{LemHar.1(2)})~Given (\ref{LemHar.1(1)(1)}), and defining $\mathcal{K}%
_{+}$ and $\Lambda$ by (\ref{eqHar.4})--(\ref{eqHar.5}), we saw that
$\mathcal{K}_{+}$ is \emph{closed,} and that, by (\ref{eqHar.2}), $\Lambda$ is
well-defined and contractive. (\ref{LemHar.1(3)})~Given (\ref{LemHar.1(1)(2)}%
), the subspace $\mathcal{K}$ in $\mathcal{H}_{0}$, defined in (\ref{eqHar.6}%
), is \emph{positive.} Indeed, if $k=k_{+}+\Lambda k_{+}$, $k_{+}%
\in\mathcal{K}_{+}$, then%
\begin{equation}
\left\langle k,Jk\right\rangle =\left\|  k_{+}\right\|  ^{2}-\left\|  \Lambda
k_{+}\right\|  ^{2}\geq0, \label{eqHar.7}%
\end{equation}
since $\Lambda$ is assumed contractive. We also easily check that
$\mathcal{K}$ in (\ref{eqHar.6}) is \emph{closed} when (\ref{LemHar.1(1)(2)})
holds, i.e., $\mathcal{K}_{+}$ is closed, and the operator $\Lambda$ in
(\ref{eqHar.3}) is contractive.
\end{proof}

\begin{corollary}
\label{CorHar.2}Let $\mathcal{K}\subset\mathcal{H}_{0}$ be a closed positive
subspace as defined in Lemma \textup{\ref{LemHar.1}} from a given $J$. Let
$\Lambda\colon\mathcal{K}_{+}\rightarrow\mathcal{H}_{-}$ be the corresponding
contraction with closed domain $\mathcal{K}_{+}\subset\mathcal{H}_{+}$, and
set%
\begin{equation}
\mathcal{N}_{+}=\left\{  k_{+}\in\mathcal{K}_{+}\mathrel{;}\Lambda^{\ast
}\Lambda k_{+}=k_{+}\right\}  . \label{eqHar.8}%
\end{equation}
Let%
\begin{equation}
\mathcal{H}_{+}\left(  \Lambda\right)  =\left(  \mathcal{K}_{+}/\mathcal{N}%
_{+}\right)  \sptilde\label{eqHar.9}%
\end{equation}
be the Hilbert space obtained by completing the quotient space $\mathcal{K}%
_{+}/\mathcal{N}_{+}$ relative to the Hilbert norm%
\begin{equation}
k_{+}\longmapsto\left\|  \left(  I-\Lambda^{\ast}\Lambda\right)  ^{\frac{1}%
{2}}k_{+}\right\|  , \label{eqHar.10}%
\end{equation}
and let%
\begin{equation}
W_{+}\colon\mathcal{K}_{+}\longrightarrow\mathcal{K}_{+}/\mathcal{N}%
_{+}\longrightarrow\mathcal{H}_{+}\left(  \Lambda\right)  \label{eqHar.11}%
\end{equation}
be the natural contractive mapping. Then%
\begin{equation}
W_{+}=P_{\mathcal{K}}P_{+}\left(  I-\Lambda^{\ast}\Lambda\right)  ^{\frac
{1}{2}}P_{+}P_{\mathcal{K}}, \label{eqHar.12}%
\end{equation}
where $P_{\mathcal{K}}$ denotes the projection of $\mathcal{H}_{0}$ onto
$\mathcal{K}$, and $P_{\pm}$ are given by \textup{(\ref{eqHar.1prime}).}
Finally there is a unitary isomorphism%
\[
T\colon\mathcal{H}_{+}\left(  \Lambda\right)  \longrightarrow\mathcal{H}%
\left(  \mathcal{K}\right)
\]
which is determined by the formula%
\begin{equation}
W=TW_{+}P_{+}P_{\mathcal{K}}. \label{eqHar.13}%
\end{equation}
\end{corollary}

\begin{proof}
Let $\mathcal{K}$ be a positive subspace, and let $\Lambda$ be the
corresponding contraction with closed domain $\mathcal{K}_{+}$, see Lemma
\ref{LemHar.1}. We saw that then $\mathcal{K}=G\left(  \Lambda\right)  $; and,
if
\begin{equation}
k=k_{+}+\Lambda k_{+},\qquad k_{+}\in\mathcal{K}_{+}, \label{eqHar.14}%
\end{equation}
then%
\begin{equation}
\left\langle k,Jk\right\rangle =\left\|  k_{+}\right\|  ^{2}-\left\|  \Lambda
k_{+}\right\|  ^{2}=\left\langle k_{+},k_{+}-\Lambda^{\ast}\Lambda
k_{+}\right\rangle =\left\|  \left(  I-\Lambda^{\ast}\Lambda\right)
^{\frac{1}{2}}k_{+}\right\|  ^{2}. \label{eqHar.15}%
\end{equation}
It follows that the assignment $k_{+}\mapsto k$ then passes to respective
quotients%
\[
\mathcal{K}_{+}/\mathcal{N}_{+}\longrightarrow\mathcal{K}/\mathcal{N},
\]
where $\mathcal{N}_{+}$ is defined in (\ref{eqHar.8}). If $T_{0}$ is the
corresponding operator $\mathcal{K}_{+}/\mathcal{N}_{+}\rightarrow
\mathcal{K}/\mathcal{N}$ induced by $k_{+}\mapsto k_{+}+\Lambda k_{+}$, then
$T_{0}$ is isometric relative to the two new norms, and it passes to the
respective completions%
\[
T=\tilde{T}_{0}\colon\underset{%
\begin{array}
[c]{c}%
\shortparallel\\
\mathcal{H}_{+}\left(  \Lambda\right)
\end{array}
}{\left(  \mathcal{K}_{+}/\mathcal{N}_{+}\right)  \sptilde}\longrightarrow
\underset{%
\begin{array}
[c]{c}%
\shortparallel\\
\mathcal{H}\left(  \mathcal{K}\right)
\end{array}
}{\left(  \mathcal{K}/\mathcal{N}\right)  \sptilde}.
\]
{}From (\ref{eqHar.14})--(\ref{eqHar.15}), we read off formula (\ref{eqHar.12}%
) for the contraction $W_{+}\colon\mathcal{K}_{+}\rightarrow\mathcal{H}%
_{+}\left(  \Lambda\right)  $. Using again (\ref{eqHar.15}), we conclude that
$T$ satisfies (\ref{eqHar.13}). Conversely, if $W$ and $W_{+}$ are constructed
from $\mathcal{K}$ and $\Lambda$, respectively, then, if we set $TW_{+}%
k_{+}=Wk$, $k\in\mathcal{K}$, then $T$ is isometric, and extends naturally to
a unitary isomorphism of $\mathcal{H}_{+}\left(  \Lambda\right)  $ onto
$\mathcal{H}\left(  \mathcal{K}\right)  $.
\end{proof}

\begin{remark}
\label{RemHar.pound}Recent work of Arveson \cite{Arv98} suggests a
multivariable version of the construction in Section \textup{\ref{Rep}} above,
i.e., reproducing kernels in several variables, as a candidate for a model in
multivariable operator theory. With this in view, one should generalize
Corollary \textup{\ref{CorHar.2}} above to the case of a system of commuting
operators $\Lambda_{i}\colon\mathcal{K}_{+}\rightarrow\mathcal{H}_{-}$,
$i=1,\dots,d$, such that%
\[
\left\|  \sum_{i=1}^{d}\Lambda_{i}k_{i}\right\|  ^{2}\leq\sum_{k=1}%
^{d}\left\|  k_{i}\right\|  ^{2}%
\]
for all $k_{1},\dots,k_{d}$, $k_{i}\in\mathcal{K}_{+}$. To make the connection
to the setup \textup{(\ref{eqHar.10})} in the present Corollary
\textup{\ref{CorHar.2},} note that the condition of Arveson is equivalent to
the operator estimate%
\[
\Lambda_{1}^{{}}\Lambda_{1}^{\ast}+\dots+\Lambda_{d}^{{}}\Lambda_{d}^{\ast
}\leq I,
\]
and the analogue of our operator from \textup{(\ref{eqHar.10})} is then
\[
\left(  I-\sum_{i=1}^{d}\Lambda_{i}^{{}}\Lambda_{i}^{\ast}\right)  ^{\frac
{1}{2}}.
\]
\end{remark}

The following observations make connections between the reflection-symmetric
operator $U$ and the subspace $\mathcal{K}$.

Let $\mathcal{H}_{+}$ and $\mathcal{H}_{-}$ be Hilbert spaces, set
\begin{equation}
\mathcal{H}_{0}=\mathcal{H}_{+}\oplus\mathcal{H}_{-},\qquad J=%
\begin{pmatrix}
I & 0\\
0 & -I
\end{pmatrix}
, \label{eqHar118.1}%
\end{equation}
and let $a\colon\mathcal{H}_{+}\rightarrow\mathcal{H}_{-}$ be an arbitrary
operator. Then set
\begin{equation}
U=U\left(  a\right)  =%
\begin{pmatrix}
a^{\ast}a & a^{\ast}\\
-a & aa^{\ast}%
\end{pmatrix}
. \label{eqHar118.2}%
\end{equation}
It follows that%
\begin{equation}
JU\left(  a\right)  J=U\left(  a\right)  ^{\ast}=U\left(  -a\right)  ,
\label{eqHar118.3}%
\end{equation}
i.e., $U\left(  a\right)  $ is reflection-symmetric. Moreover, $U=U\left(
a\right)  $ satisfies%
\[
U^{\ast}U=%
\begin{pmatrix}
a^{\ast}a + \left(  a^{\ast}a \right)  ^{2} & 0\\
0 & aa^{\ast}+\left(  aa^{\ast}\right)  ^{2}%
\end{pmatrix}
.
\]
Conversely, every operator $U\colon\mathcal{H}_{0}\rightarrow\mathcal{H}_{0}$
which satisfies%
\begin{equation}
JUJ=U^{\ast}, \label{eqHar118.4}%
\end{equation}
and%
\begin{equation}
U^{\ast}U=\left(
\begin{tabular}
[c]{c|c}%
$\operatorname*{operator}_{1}$ & $0$\\\hline
$0$ & $\operatorname*{operator}_{2}$%
\end{tabular}
\right)  \label{eqHar118.5}%
\end{equation}
relative to the decomposition (\ref{eqHar118.1}) is of the form
\begin{equation}
U=%
\begin{pmatrix}
s_{1} & a^{\ast}\\
-a & s_{2}%
\end{pmatrix}
\label{eqHar118.6}%
\end{equation}
for some operator $a\colon\mathcal{H}_{+}\rightarrow\mathcal{H}_{-}$, and for
two selfadjoint operators $s_{1}$ and $s_{2}$ in the respective Hilbert spaces
$\mathcal{H}_{+}$ and $\mathcal{H}_{-}$, and satisfying the intertwining
relation:%
\begin{equation}
as_{1}=s_{2}a. \label{eqHar118.7}%
\end{equation}

Returning to the classical example from Section \ref{Int} above, let
$\mathcal{H}_{0}:=L^{2}\left(  \mathbb{T}\right)  $, and set%
\begin{equation}
Jf\left(  z\right)  :=f\left(  \bar{z}\right)  ,\qquad f\in L^{2}\left(
\mathbb{T}\right)  ,\;z\in\mathbb{T}. \label{eqHar.16}%
\end{equation}

\begin{proposition}
\label{ProHar.query}Let $H^{2}=H^{2}\left(  \mathbb{T}\right)  $, and
$H^{\infty}=H^{\infty}\left(  \mathbb{T}\right)  $ be the usual Hardy spaces
of harmonic analysis. Let $b\in H^{\infty}$ be given, and suppose that
$\left\|  b\right\|  _{\infty}\leq1$. Define the subspace $\mathcal{K}\left(
b\right)  \subset\mathcal{H}_{0}$ \textup{(}$=L^{2}\left(  \mathbb{T}\right)
$\textup{)} as follows:%
\begin{equation}
\mathcal{K}\left(  b\right)  =\left\{  \left(  1-b\left(  \bar{z}\right)
\right)  k\left(  \bar{z}\right)  +\left(  1+b\left(  z\right)  \right)
k\left(  z\right)  \mathrel{;}k\in H^{2}\right\}  \label{eqHar.17}%
\end{equation}
Then $\mathcal{K}\left(  b\right)  $ is a maximal positive subspace of
$\mathcal{H}_{0}$ relative to the given reflection operator $J$ from
\textup{(\ref{eqHar.16}).} Moreover, the space $\mathcal{K}\left(  b\right)  $
is invariant under the shift%
\begin{equation}
Uf\left(  z\right)  =zf\left(  z\right)  ,\qquad f\in L^{2}\left(
\mathbb{T}\right)  ,\;z\in\mathbb{T},\label{eqHar.18}%
\end{equation}
if and only if $b\equiv1$. In that case, $\mathcal{H}\left(  \mathcal{K}%
\right)  $ is one-dimensional, and $S\left(  U\right)  =0$.
\end{proposition}

\begin{proof}
The proof is based on Corollary \ref{CorHar.2} above. Since $J$ is given by
(\ref{eqHar.16}) at the outset, the two subspaces $\mathcal{H}_{\pm}\subset
L^{2}\left(  \mathbb{T}\right)  $ are then determined from (\ref{eqHar.1prime}%
), applied to $J$. Let $\mathcal{K}=H^{2}\left(  \mathbb{T}\right)  $, and set
$\mathcal{K}_{\pm}:=P_{\pm}\mathcal{K}$. Then $\mathcal{K}_{\pm}%
=\mathcal{H}_{\pm}$, where
\begin{equation}
\mathcal{K}_{\pm}=\left\{  k\left(  z\right)  \pm k\left(  \bar{z}\right)
\mathrel{;}k\in H^{2}\right\}  . \label{eqHar.19}%
\end{equation}
Let $b\in H^{\infty}$, $\left\|  b\right\|  _{\infty}\leq1$, be given, and
define $\Lambda=\Lambda_{b}$ by%
\begin{equation}
\Lambda\left(  P_{+}k\right)  :=P_{-}\left(  bk\right)  \text{,\qquad for all
}k\in H^{2}. \label{eqHar.20}%
\end{equation}
Then it follows from $\mathcal{K}_{+}=\mathcal{H}_{+}$ that $\Lambda$ is a
contractive operator with domain $\mathcal{H}_{+}$ and mapping into
$\mathcal{H}_{-}$. The corresponding positive subspace, see Lemma
\ref{LemHar.1}, is that which is given by (\ref{eqHar.17}). The space
$\mathcal{K}\left(  b\right)  $ is maximally positive. A positive subspace
$\mathcal{K}^{\prime}$ satisfying $\mathcal{K}\left(  b\right)  \subset
\mathcal{K}^{\prime}$ would correspond to a contractive operator
$\Lambda^{\prime}$ mapping $\mathcal{H}_{+}$ into $\mathcal{H}_{-}$ and
extending $\Lambda$, in the sense that the graph of $\Lambda^{\prime}$
contains that of $\Lambda$. But then $\Lambda=\Lambda^{\prime}$ and therefore
$\mathcal{K}\left(  b\right)  =\mathcal{K}^{\prime}$ by the uniqueness part in
Lemma \ref{LemHar.1}. This proves that $\mathcal{K}\left(  b\right)  $ is
maximally positive in $L^{2}\left(  \mathbb{T}\right)  $.

The contractive property for the operator $\Lambda=\Lambda_{b}$ in
(\ref{eqHar.20}) follows from the two assumptions on $b$, i.e., $b\in
H^{\infty}$, and $\left\|  b\right\|  _{\infty}\leq1$. Indeed, if $k\in H^{2}%
$, then%
\begin{multline*}
\left\|  P_{-}\left(  bk\right)  \right\|  _{2}^{2}=\left\|  \frac{1}%
{2}\left(  -b\left(  \bar{z}\right)  k\left(  \bar{z}\right)  +b\left(
z\right)  k\left(  z\right)  \right)  \right\|  _{2}^{2}\\
=\frac{1}{2}\left(  \left\|  bk\right\|  _{2}^{2}-\left|  b\left(  0\right)
k\left(  0\right)  \right|  ^{2}\right) \\
\leq\frac{1}{2}\left\|  bk\right\|  _{2}^{2}\leq\frac{1}{2}\left\|  b\right\|
_{\infty}^{2}\left\|  k\right\|  _{2}^{2}\leq\frac{1}{2}\left\|  k\right\|
_{2}^{2}=\left\|  P_{+}k\right\|  _{2}^{2}.
\end{multline*}
This proves that the operator $\Lambda=\Lambda_{b}$ in (\ref{eqHar.20}) is
indeed well-defined and contractive. We then conclude from Lemma
\ref{LemHar.1}(\ref{LemHar.1(3)}) that the corresponding positive subspace
$\mathcal{K}\left(  b\right)  $ is the graph of $\Lambda_{b}$. An application
of (\ref{eqHar.6}) from Lemma \ref{LemHar.1} then finally yields
(\ref{eqHar.17}) as claimed.

If it were the case that $\mathcal{K}\left(  b\right)  $ ($=G\left(
\Lambda_{b}\right)  $) were invariant under the shift $U$ of (\ref{eqHar.18}),
then from Beurling's theorem, there would be a unitary function $u\in
L^{\infty}\left(  \mathbb{T}\right)  $ such that%
\begin{equation}
\mathcal{K}\left(  b\right)  =uH^{2}. \label{eqHar.21}%
\end{equation}
(Recall $u\in L^{\infty}$ is said to be unitary if the corresponding
multiplication operator $M_{u}$ on $L^{2}$ is unitary.) But identity in
(\ref{eqHar.21}) for some unitary $u\in L^{\infty}$ is possible only if the
factor $\left(  1-b\left(  \bar{z}\right)  \right)  $ in (\ref{eqHar.17})
vanishes identically on $\mathbb{T}$, and it follows therefore that
$\mathcal{K}\left(  b\right)  $ can only be shift-invariant if $b\equiv1$. In
this case, $\mathcal{K}\left(  b\right)  =\mathcal{K=}H^{2}$ reduces to the
special case which we studied in Section \ref{Int}. In that case, the
contraction $\Lambda$ from (\ref{eqHar.20}) reduces to $\Lambda\left(
P_{+}k\right)  =P_{-}k$, and%
\[
\left\langle k,Jk\right\rangle =\left\|  P_{+}k\right\|  ^{2}-\left\|
P_{-}k\right\|  ^{2}=\left|  c_{0}\right|  ^{2}\text{\quad if\quad}k\left(
z\right)  =\sum_{n=0}^{\infty}c_{n}z^{n}\in H^{2}.
\]
Hence $\mathcal{H}\left(  \mathcal{K}\right)  $ is one-dimensional. Since%
\[
Uk\left(  z\right)  =zk\left(  z\right)  =c_{0}z+c_{1}z^{2}+\cdots
\]
has zero constant term, the selfadjoint operator $S\left(  U\right)  $ on
$\mathcal{H}\left(  \mathcal{K}\right)  $, induced from $U$, is zero, and the
proof is completed.
\end{proof}

Elaborating on the abstract setup in Proposition \ref{ProInt.3}, we conclude
with a family of finite-dimensional positive subspaces in $H^{2}\oplus H^{2}$.

The simplest situation when a triple $\left(  \mathcal{H}_{0},\mathcal{K}%
,J\right)  $ arises in an application is the case of the Pick--Nevanlinna
interpolation problem. In that case, let%
\[
\mathcal{H}_{0}=\ell_{+}^{2}\oplus\ell_{+}^{2},\qquad J=%
\begin{pmatrix}
I & 0\\
0 & -I
\end{pmatrix}
,
\]
$N\in\mathbb{N}$, distinct points $z_{1},\dots,z_{N}\in D=\left\{
z\in\mathbb{C}\mathrel{;}\left|  z\right|  <1\right\}  $, and $w_{1}%
\dots,w_{N}\in\mathbb{C}$, be given. The Pick--Nevanlinna theorem states that
there exists a function $\varphi\in H^{\infty}\left(  D\right)  $ such that
$\varphi\left(  z_{i}\right)  =w_{i}$ for each $i$, and $\left\|
\varphi\right\|  _{\infty}\leq$ $1$ if and only if the corresponding $N\times
N$ matrix $\left(  \frac{1-\bar{w}_{i}w_{j}}{1-\bar{z}_{i}z_{j}}\right)  $ is
positive semidefinite. We will now assume the latter, and relate it to the
$\mathcal{K}$-problem. Then set%
\[
\mathcal{K}:=\left\{
\begin{pmatrix}
\left(  \sum_{i}c_{i}^{{}}z_{i}^{n}\right)  _{n=0_{\mathstrut}}^{\infty}\\
\left(  \sum_{i}c_{i}^{{}}w_{i}^{{}}z_{i}^{n}\right)  _{n=0}^{\infty
^{\mathstrut}}%
\end{pmatrix}
\mathrel{;}c_{1},c_{2},\dots,c_{N}\in\mathbb{C}\right\}  \subset%
\begin{pmatrix}
\ell_{+}^{2}\\
\ell_{+}^{2}%
\end{pmatrix}
^{\oplus}.
\]
It is an $N$-dimensional subspace, and so closed. For general vectors
$k=k\left(  c\right)  $, $c=\left(  c_{1},\dots,c_{N}\right)  $ in
$\mathcal{K}$, the term $\left\langle k,Jk\right\rangle =\left\|
P_{+}k\right\|  ^{2}-\left\|  P_{-}k\right\|  ^{2}$ computes out as%
\[
\sum_{n}\left|  \sum_{i}c_{i}^{{}}z_{i}^{n}\right|  ^{2}-\sum_{n}\left|
\sum_{i}c_{i}^{{}}w_{i}^{{}}z_{i}^{n}\right|  ^{2}=\sum_{i}\sum_{j}%
\frac{1-\bar{w}_{i}w_{j}}{1-\bar{z}_{i}z_{j}}\bar{c}_{i}c_{j}\geq0,
\]
assuming the Pick--Nevanlinna condition.

Since we also work with the $H^{2}$-version of $\ell_{+}^{2}$, we note that
the above positive subspace $\mathcal{K}$ has an equivalent form in
$\mathcal{H}_{0}=H^{2}\oplus H^{2}$. There we have the reproducing kernel
$q_{z}\left(  \zeta\right)  =\left(  1-\bar{z}\zeta\right)  ^{-1}$, and
$\mathcal{K}$ then takes the form of column vectors as follows:%
\[
\mathcal{K}=\left\{
\begin{pmatrix}
\sum_{i}c_{i}q_{z_{i}}\\
\sum_{i}c_{i}w_{i}q_{z_{i}}%
\end{pmatrix}
\mathrel{;}c_{1},c_{2},\dots,c_{N}\in\mathbb{C}\right\}  .
\]

The Pick--Nevanlinna problem was stated in terms of the pair $\mathcal{K}$,
$J=\left(  I\oplus\left(  -I\right)  \right)  $, but if we use instead
$J=\left(
\begin{smallmatrix}
0 & I\\
I & 0
\end{smallmatrix}
\right)  $, then it is easy to check that the corresponding condition,
$\left\langle k,Jk\right\rangle \geq0$ for $k\in\mathcal{K}$, is now
equivalent to the matrix order relation, $\left(  \frac{\bar{w}_{i}+w_{j}%
}{1-\bar{z}_{i}z_{j}}\right)  \geq0$, i.e., equivalent to%
\[
\sum_{i=1}^{N}\sum_{j=1}^{N}\bar{c}_{i}\left(  \frac{\bar{w}_{i}+w_{j}}%
{1-\bar{z}_{i}z_{j}}\right)  c_{j}\geq0\text{\qquad for all }c_{1},\dots
,c_{N}\in\mathbb{C}.
\]
This alternative is in turn equivalent to a solution to the interpolation
problem $\varphi\left(  z_{i}\right)  =w_{i}$ for each $i$, and
$\operatorname{Re}\varphi\geq0$ in $D$ for some interpolating analytic
function $\varphi$. Hence both of the classical interpolation problems
correspond to positivity for a pair $\left(  \mathcal{K},J\right)  $ where
$\mathcal{K}\subset H^{2}\oplus H^{2}$ is as stated, but where $J$ changes
from one problem to the other.

A nice solution to both problems is presented in the classic paper
\cite{Sar67}. (See also \cite{FaKo94}.)

\section{\label{Han}Hankel operators}

In this section, we consider the direct sum of the unilateral shift $A$ and
its adjoint $A^{\ast}$, i.e., $U=A\oplus A^{\ast}$. If $J=\left(
\begin{smallmatrix}
0 & I\\
I & 0
\end{smallmatrix}
\right)  $, then $JUJ=U^{\ast}$, and we solve the problem of finding the
subspaces $\mathcal{K}\subset\ell_{+}^{2}\oplus\ell_{+}^{2}$ which satisfy the
positivity (\ref{eqThmRef.1(2)}) of Theorem \ref{ThmRef.1}, and are invariant
under $U$. This is analogous to (and yet very different from) the classical
solution of Beurling \cite[chapter 6]{Hel95} which gives the invariant
subspaces for $A$. Recall the invariant subspaces for $A$ are in $1$--$1$
correspondence with the \emph{inner functions,} i.e., functions $\xi\in
H^{\infty}$ such that $\left|  \xi\left(  e^{i\theta}\right)  \right|  =1$
$\mathrm{a.e.}$ $\theta\in\left[  -\pi,\pi\right)  $. For our present problem
with $A\oplus A^{\ast}$, we will first reduce the analysis to considering
closed invariant subspaces $\mathcal{K}\subset\ell_{+}^{2}\oplus\ell_{+}^{2}$
which are maximally positive. This reduction follows in fact from an
application of Beurling's theorem. We then show that those invariant subspaces
$\mathcal{K}$ are in $1$--$1$ correspondence with positive and finite Borel
measures $\mu$ on $\left[  -1,1\right]  $ in such a way that the corresponding
induced selfadjoint operator $S_{\mu}\left(  A\oplus A^{\ast}\right)  $,
acting on $\mathcal{H}\left(  \mathcal{K}\right)  $, is unitarily equivalent
to multiplication by the real variable $x$ on $L_{\mu}^{2}\left(  \left[
-1,1\right]  \right)  $, i.e., $f\left(  x\right)  \mapsto xf\left(  x\right)
$, on the $L^{2}$ space given by $\int_{-1}^{1}\left|  f\left(  x\right)
\right|  ^{2}\,d\mu\left(  x\right)  <\infty$, and defined from a finite
positive measure $\mu$ on $\left[  -1,1\right]  $. We also make explicit how a
subspace $\mathcal{K}=\mathcal{K}_{\mu}$ with the desired properties may be
reconstructed from some given measure $\mu$ as specified.

We first give some Hilbert-space background: Let $\mathcal{H}$ be a Hilbert
space, and let $A$ be a bounded operator in $\mathcal{H}$. Then%
\begin{equation}
U:=%
\begin{pmatrix}
A & 0\\
0 & A^{\ast}%
\end{pmatrix}
\text{\qquad on }\mathcal{H}_{0}:=\mathcal{H}\oplus\mathcal{H} \label{eqHan.1}%
\end{equation}
satisfies%
\begin{equation}
JUJ=U^{\ast} \label{eqHan.2}%
\end{equation}
relative to%
\begin{equation}
J=%
\begin{pmatrix}
0 & I\\
I & 0
\end{pmatrix}
, \label{eqHan.3}%
\end{equation}
i.e., the operator $J$ on $\mathcal{H}_{0}$ is given by $J\left(  h\oplus
k\right)  =k\oplus h$. This observation also shows that the identity
(\ref{eqHan.2}) typically does not imply any special property for the
operators making up $U$. On the other hand, the example in Section \ref{Ref}
had $U$ unitary relative to the original Hilbert space $\mathcal{H}_{0}$.

We wish to compute the correspondence $U\mapsto S\left(  U\right)  $ of
Theorem \ref{ThmRef.1} in the case of (\ref{eqHan.1}) and (\ref{eqHan.3}).
Given a subspace $\mathcal{K}\subset\mathcal{H}_{0}$ such that%
\begin{equation}
U\left(  \mathcal{K}\right)  \subset\mathcal{K}, \label{eqHan.3prime}%
\end{equation}
we will pass to the new Hilbert space
\begin{equation}
\mathcal{H}\left(  \mathcal{K}\right)  =\left(  \mathcal{K}/\mathcal{N}%
\right)  \sptilde, \label{eqHan.3primeprime}%
\end{equation}
where $\mathcal{N}=\left\{  k\in\mathcal{K}\mathrel{;}\left\langle
k,Jk\right\rangle =0\right\}  $. We say that $\mathcal{K}$ is the \emph{graph}
of some operator from a domain $D\left(  \Gamma\right)  \subset\mathcal{H}$
into $\mathcal{H}$, if%
\begin{equation}%
\begin{pmatrix}
0\\
h
\end{pmatrix}
\in\mathcal{K}\Longrightarrow h=0. \label{eqHan.3primeprimeprime}%
\end{equation}
But in view of (\ref{eqHan.3}), vectors of the form $\left(
\begin{smallmatrix}
0\\
h
\end{smallmatrix}
\right)  $ are automatically in $\mathcal{N}$, and so do not contribute to
$\mathcal{H}\left(  \mathcal{K}\right)  $ of (\ref{eqHan.3primeprime}). We
will suppose, therefore, that the spaces $\mathcal{K}$ of (\ref{eqHan.3prime})
have the form $\mathcal{K}=G\left(  \Gamma\right)  $. Note that the operator
$\Gamma$ of which $\mathcal{K}$ is the graph need not have dense domain. The
subspace $\mathcal{K}$ is said to be \emph{positive} if $\left\langle
k,Jk\right\rangle \geq0$ for $k\in\mathcal{K}$, and \emph{maximally positive}
if it is maximal (relative to inclusion) with respect to this property. It
follows from (\ref{eqHan.3}) that the maximally positive subspaces
$\mathcal{K}$ of the form $\mathcal{K}=G\left(  \Gamma\right)  $ correspond to
operators $\Gamma$ which are \emph{dissipative, closed,} and have \emph{dense
domain} in $\mathcal{H}$. The corresponding Cayley transform%
\begin{equation}
\Lambda:=\left(  I-\Gamma\right)  \left(  I+\Gamma\right)  ^{-1}
\label{eqHan.3iv}%
\end{equation}
is then contractive and everywhere defined on $\mathcal{H}$, and it
corresponds to the contraction also denoted $\Lambda$ from Lemma
\ref{LemRef.2}. This contraction derives from the general contractive
transformation%
\begin{equation}
P_{+}k\longmapsto P_{-}k,\qquad k\in\mathcal{K}, \label{eqHan.3v}%
\end{equation}
where $P_{\pm}=\frac{1}{2}\left(  I\pm J\right)  $. Using (\ref{eqHan.3}) we
get%
\[
P_{\pm}%
\begin{pmatrix}
h\\
\Gamma h
\end{pmatrix}
=\frac{1}{2}%
\begin{pmatrix}
h\pm\Gamma h\\
h\pm\Gamma h
\end{pmatrix}
\text{\qquad for }h\in D\left(  \Gamma\right)  ,
\]
and so%
\[
\left\|  P_{\pm}%
\begin{pmatrix}
h\\
\Gamma h
\end{pmatrix}
\right\|  =\frac{1}{\sqrt{2}}\left\|  h\pm\Gamma h\right\|  .
\]
Since (\ref{eqHan.3v}) is contractive, it follows that $\Lambda$ in
(\ref{eqHan.3iv}) is well-defined and also contractive. Let $A$ in
(\ref{eqHan.1}) be the unilateral shift. Then of course $U$ will not even be
normal. Nonetheless, the possibilities for reflection symmetry yield a richer
family, and we will show here that the possibilities can even be classified,
i.e., if $A$ in (\ref{eqHan.1}) is the unilateral shift.

Let $\mathcal{H}=H^{2}$. We will use both of the representations $f\left(
z\right)  =\sum_{n=0}^{\infty}c_{n}z^{n}$, and $\left(  c_{0,}c_{1}%
,c_{2},\dots\right)  $ for elements in $H^{2}$, i.e., the function vs.\ its
Fourier series. Hence $A$ takes alternately the form%
\begin{equation}
\left(  Af\right)  \left(  z\right)  =zf\left(  z\right)  ,\qquad f\in
H^{2},\;z\in\mathbb{T}, \label{eqHan.4}%
\end{equation}
or%
\begin{equation}
A\left(  c_{0},c_{1},c_{2},\dots\right)  =\left(  0,c_{0},c_{1},c_{2}%
,\dots\right)  ,\qquad\left(  c_{n}\right)  _{n=0}^{\infty}\in\ell^{2},
\label{eqHan.5}%
\end{equation}
and $A^{\ast}$ given by $A^{\ast}\left(  c_{0},c_{1},c_{2},\dots\right)
=\left(  c_{1},c_{2},c_{3},\dots\right)  $.

It is immediate that, if $\Gamma$ is an operator in $\mathcal{H}=H^{2}$, with
domain $D\left(  \Gamma\right)  $, and graph $G\left(  \Gamma\right)
=\left\{  \left(
\begin{smallmatrix}
h\\
\Gamma h
\end{smallmatrix}
\right)  \mathrel{;}h\in D\left(  \Gamma\right)  \right\}  $, then
$\mathcal{K}:=G\left(  \Gamma\right)  $ satisfies the positivity%
\begin{equation}
\left\langle k,Jk\right\rangle \geq0\text{\qquad for all }k\in\mathcal{K}
\label{eqHan.6}%
\end{equation}
if and only if $\Gamma$ is \emph{dissipative,} meaning
\begin{equation}
\operatorname{Re}\left\langle h,\Gamma h\right\rangle \geq0\text{\qquad for
all }h\in D\left(  \Gamma\right)  . \label{eqHan.7}%
\end{equation}

It is easy to show, see, e.g., \cite{Phil}, that if $\Gamma$ is dissipative,
then the closure of $G\left(  \Gamma\right)  $, i.e., $\overline{G\left(
\Gamma\right)  }$, is also the graph of a dissipative operator, denoted
$\bar{\Gamma}$. (An operator is said to be \emph{closed} if its graph is
closed.) We will consider subspaces $\mathcal{K}$ which are invariant under
$U=\left(
\begin{smallmatrix}
A & 0\\
0 & A^{\ast}%
\end{smallmatrix}
\right)  $. But if $\mathcal{K}$ is invariant, then so is $\overline
{\mathcal{K}}$, and we will restrict attention to closed subspaces, and
corresponding closed operators.

\begin{lemma}
\label{LemHan.1}Let $U=\left(
\begin{smallmatrix}
A & 0\\
0 & A^{\ast}%
\end{smallmatrix}
\right)  $ be built from the shift $A$, see \textup{(\ref{eqHan.4}),} and let
$\Gamma$ be an operator with domain $D\left(  \Gamma\right)  $ in $H^{2}$, and
graph $G\left(  \Gamma\right)  $ in $H^{2}\oplus H^{2}$. Then%
\begin{equation}
U\left(  G\left(  \Gamma\right)  \right)  \subset G\left(  \Gamma\right)
\label{eqHan.8}%
\end{equation}
if and only if $D\left(  \Gamma\right)  $ is $A$-invariant and
\begin{equation}
\Gamma A=A^{\ast}\Gamma\text{\qquad on }D\left(  \Gamma\right)  .
\label{eqHan.9}%
\end{equation}
\end{lemma}

\begin{proof}
Since%
\[
U%
\begin{pmatrix}
h\\
\Gamma h
\end{pmatrix}
=%
\begin{pmatrix}
Ah\\
A^{\ast}\Gamma h
\end{pmatrix}
\text{\qquad for }h\in D\left(  \Gamma\right)  ,
\]
we see that (\ref{eqHan.8}) holds if and only if $\Gamma Ah=A^{\ast}\Gamma h$,
which is the conclusion.
\end{proof}

However, the operators $\Gamma$ satisfying (\ref{eqHan.9}) are the Hankel
operators. Relative to the standard basis in $H^{2}$, such a $\Gamma$ has the
form%
\begin{equation}
\left(  \Gamma x\right)  _{n}=\sum_{m=0}^{\infty}\gamma_{n+m}x_{m}
\label{eqHan.10}%
\end{equation}
for $n=0,1,\dots$, where $\gamma$ is some sequence, $\gamma\in\ell^{2}$. While
the bounded Hankel operators are known, the interesting ones, for reflection
positivity, will be unbounded ones. (Recall $\Gamma=\Gamma_{\gamma}$ is
bounded in $H^{2}\ $if and only if there is some $\varphi\in L^{\infty}\left(
\mathbb{T}\right)  $ such that $\gamma_{n}=\hat{\varphi}\left(  -n\right)  $,
$n=0,1,\dots$, see, e.g., \cite{Pow82}.)

While we can reduce to the case when $\mathcal{K}=G\left(  \Gamma\right)  $ is
closed in $H^{2}\oplus H^{2}$, the domain $D\left(  \Gamma\right)  $ is not
closed in $H^{2}$, but only dense.

\begin{lemma}
\label{LemHan.2}Let $\Gamma=\Gamma_{\gamma}$ be the closed operator defined in
\textup{(\ref{eqHan.10})} when it is assumed that%
\begin{equation}
\operatorname{Re}\gamma_{n}\geq0\text{\qquad for all }n=0,1,2,\dots.
\label{eqHan.12}%
\end{equation}
Then%
\begin{equation}
\left(  I+\Gamma\right)  D\left(  \Gamma\right)  =H^{2}. \label{eqHan.13}%
\end{equation}
\end{lemma}

\begin{proof}
It follows from (\ref{eqHan.10}) that the condition (\ref{eqHan.12}) on the
sequence $\left(  \gamma_{n}\right)  _{n=0}^{\infty}$ is equivalent to
$\Gamma_{\gamma}$ being dissipative. Hence, since $\left(  \gamma_{n}\right)
\in\ell^{2}$, the operator $\Gamma$ has a dense domain $D\left(
\Gamma\right)  $, and the closure of $\Gamma$ is well-defined. We will work
with the closure, and refer to $\Gamma$ as the closed operator. Notice that if
$\Gamma$ is defined from a sequence $\left(  \gamma_{n}\right)  $, then the
adjoint operator $\Gamma^{\ast}$ is defined from the sequence $\left(
\bar{\gamma}_{n}\right)  $; and so, by (\ref{eqHan.12}), both are dissipative.
In particular,%
\begin{equation}
\operatorname{Re}\left\langle h,\Gamma^{\ast}h\right\rangle \geq0
\label{eqHan.14}%
\end{equation}
for all $h\in D\left(  \Gamma^{\ast}\right)  $. To prove (\ref{eqHan.13}),
suppose $h\perp\left(  I+\Gamma\right)  D\left(  \Gamma\right)  $. Then $h\in
D\left(  \Gamma^{\ast}\right)  $, and $\Gamma^{\ast}h=-h$. Since then
$\operatorname{Re}\left\langle h,\Gamma^{\ast}h\right\rangle =-\left\|
h\right\|  ^{2}$, this contradicts (\ref{eqHan.14}), unless $h=0$. Hence
$\left(  I+\Gamma\right)  D\left(  \Gamma\right)  $ is dense in $H^{2}$. But
it is also closed since $\Gamma$ is closed and dissipative.
\end{proof}

\begin{theorem}
\label{ThmHan.3}The maximally positive subspaces $\mathcal{K}\subset
H^{2}\oplus H^{2}$ which are invariant under $U=\left(
\begin{smallmatrix}
A & 0\\
0 & A^{\ast}%
\end{smallmatrix}
\right)  $, $A$ the unilateral shift, have the form%
\begin{equation}
\mathcal{K}=G\left(  \Gamma_{\gamma}\right)  \mod{\mathcal{N}},
\label{eqHan.15}%
\end{equation}
where the sequence $\gamma\in\ell^{2}$ satisfies
\begin{equation}
2\operatorname{Re}\gamma_{n}=\int_{\mathbb{R}}x^{n}\,d\mu\left(  x\right)
\label{eqHan.16}%
\end{equation}
for some positive and finite Borel measure $\mu$ on the interval $I=\left[
-1,1\right]  \subset\mathbb{R}$. If $\mathcal{K}$ comes from such a measure
$\mu$, then $\mu$ is unique, and the pair $\left(  \mathcal{H}\left(
\mathcal{K}\right)  ,S\left(  U\right)  \right)  $ may be taken to be
$L^{2}\left(  I,d\mu\right)  $ for the Hilbert space $\mathcal{H}\left(
\mathcal{K}_{\mu}\right)  $, and multiplication by $x$ on $L^{2}\left(
I,d\mu\right)  $ for the induced selfadjoint operator $S_{\mu}\left(
U\right)  $.
\end{theorem}

\begin{proof}
We begin with a lemma.

\begin{lemma}
\label{LemHan154}Let $\gamma_{n}\in\mathbb{C}$, $n=0,1,2,\dots$, be a sequence
such that all the sums%
\[
S_{\gamma}\left(  \zeta\right)  :=\sum_{n}\sum_{m}\bar{\zeta}_{n}\gamma
_{n+m}\zeta_{m}%
\]
satisfy $S_{\gamma}\left(  \zeta\right)  \geq0$ for sequences $\left(
\zeta_{n}\right)  $ which are eventually zero. Let $\mu$ be a positive Borel
measure on $I:=\left[  -1,1\right]  $ with finite moments%
\[
\gamma_{n}=\int_{-1}^{1}x^{n}\,d\mu\left(  x\right)  ,\qquad n=0,1,2,\dots.
\]
Let $\Gamma$ be the \textup{(}possibly unbounded\/\textup{)} Hankel operator
with symbol sequence $\left(  \gamma_{n}\right)  $.

\begin{enumerate}
\item \label{LemHan154(1)}Then the following are equivalent:

\begin{enumerate}
\item \label{LemHan154(1)(1)}$\openone\in D\left(  \Gamma\right)  $,

\item \label{LemHan154(1)(2)}$e_{n}\left(  z\right)  :=z^{n}\in D\left(
\Gamma\right)  $ for \emph{some} $n\in\left\{  0,1,2,\dots\right\}  $,

\item \label{LemHan154(1)(3)}$e_{n}\left(  z\right)  :=z^{n}\in D\left(
\Gamma\right)  $ for \emph{all} $n\in\left\{  0,1,2,\dots\right\}  $, and

\item \label{LemHan154(1)(4)}$\left(  \gamma_{n}\right)  _{n=0}^{\infty}%
\in\ell^{2}$.
\end{enumerate}

\item \label{LemHan154(2)}If one, and therefore all, the conditions hold, then%
\[
\lim_{n\rightarrow\infty}\left\|  \Gamma\left(  e_{n}\right)  \right\|  =0.
\]

\item \label{LemHan154(3)}The conditions are satisfied if%
\begin{equation}
\int_{-1}^{1}\left(  1-x^{2}\right)  ^{-\frac{1}{2}}\,d\mu\left(  x\right)
<\infty. \label{eqHan154.1}%
\end{equation}
But \textup{(\ref{eqHan154.1})} is more restrictive than
\textup{(\ref{LemHan154(1)(1)})--(\ref{LemHan154(1)(4)})} in
\textup{(\ref{LemHan154(1)}).}
\end{enumerate}
\end{lemma}

\begin{proof}
We view $\Gamma=\Gamma_{\gamma}$ as an operator on $H^{2}$, and note that, if
$z^{n}\in D\left(  \Gamma\right)  $, then%
\[
\Gamma\left(  z^{n}\right)  =\sum_{m=0}^{\infty}\gamma_{n+m}z^{m}.
\]
Equivalently, setting $e_{n}\left(  z\right)  :=z^{n}$,
\[
\Gamma\left(  e_{n}\right)  \left(  z\right)  =\sum_{m}\gamma_{n+m}z^{m}.
\]
The equivalence of conditions (\ref{LemHan154(1)(1)})--(\ref{LemHan154(1)(4)})
of (\ref{LemHan154(1)}) is immediate from this. Indeed, if $e_{n}\in D\left(
\Gamma\right)  $, then $\left\|  \Gamma\left(  e_{n}\right)  \right\|
^{2}=\sum_{k=n}^{\infty}\left|  \gamma_{k}\right|  ^{2}$. So this decides
(\ref{LemHan154(1)(4)}); and (\ref{LemHan154(2)}) also follows. Hence for
(\ref{LemHan154(3)}), it is enough to show that (\ref{LemHan154(1)(1)})
follows from (\ref{eqHan154.1}). Let $\left(  c_{0},c_{1},\dots\right)  $ be a
sequence which is eventually zero. Then
\begin{multline*}
\left|  \sum_{n=0}^{\infty}\gamma_{n}c_{n}\right|  =\left|  \sum_{n=0}%
^{\infty}\int_{I}x^{n}c_{n}\,d\mu\left(  x\right)  \right|  \leq\int_{I}%
\sum_{n=0}^{\infty}\left|  x^{n}c_{n}\right|  \,d\mu\left(  x\right) \\
\leq\int_{I}\left(  \sum_{n=0}^{\infty}x^{2n}\right)  ^{\frac{1}{2}}\left(
\sum_{n=0}^{\infty}\left|  c_{n}\right|  ^{2}\right)  ^{\frac{1}{2}}%
\,d\mu\left(  x\right)  =\left\|  \left(  c_{n}\right)  \right\|  _{\ell^{2}%
}\cdot\int_{I}\left(  1-x^{2}\right)  ^{-\frac{1}{2}}\,d\mu\left(  x\right)  ,
\end{multline*}
and the integral on the right is finite by assumption (\ref{eqHan154.1}). It
follows that the sequence $\left(  \gamma_{n}\right)  $ defines a bounded
linear functional on $H^{2}\simeq\ell_{+}^{2}$, and so it is in $\ell_{+}^{2}$
by Riesz's theorem. Equivalently, $\Gamma\left(  e_{0}\right)  \left(
z\right)  =\sum_{n=0}^{\infty}\gamma_{n}z^{n}$ defines an element of $H^{2}$,
and so (\ref{LemHan154(1)(1)}) holds, and in fact $\Gamma_{\gamma}$ is densely
defined as an operator on $H^{2}$.
\end{proof}

We now continue with the proof of Theorem \ref{ThmHan.3}. Let $\mathcal{K}$ be
given, and assume it has the properties stated in the theorem. Then from
Theorem \ref{ThmRef.1}, we know that there is a selfadjoint version $S\left(
U\right)  $ in a Hilbert space $\mathcal{H}\left(  \mathcal{K}\right)  $. With
the data from Theorem \ref{ThmRef.1}, we also know that the pair $\left(
\mathcal{H}\left(  \mathcal{K}\right)  ,S\left(  U\right)  \right)  $ is
unique up to unitary equivalence. Since the spectral radius of $U$ in the
present theorem is clearly one, we get, from Theorem \ref{ThmRef.1}%
(\ref{eqThmRef.1(5)}), that $\left\|  S\left(  U\right)  \right\|  \leq1$.
Suppose for the moment that $S\left(  U\right)  $ is realized as
multiplication by $x$ on $L^{2}\left(  \mathbb{R},\mu\right)  $. Then the
spectrum of $S_{\mu}\left(  U\right)  $ must be contained in $I=\left[
-1,1\right]  $, and so the support of $\mu$ must be contained in $I$.

We saw in Lemmas \ref{LemHan.1} and \ref{LemHan.2} that $\mathcal{K}$ must
have the desired form (\ref{eqHan.15}) for some dissipative operator $\Gamma$
with dense domain $D\left(  \Gamma\right)  $ in $\mathcal{H}$. Since $G\left(
\Gamma\right)  $ is mapped into itself by $\left(
\begin{smallmatrix}
A & 0\\
0 & A^{\ast}%
\end{smallmatrix}
\right)  $, we get the commutation identity (\ref{eqHan.9}). Writing out the
positivity (\ref{eqHan.6}) for $k=\left(
\begin{smallmatrix}
h\\
\Gamma h
\end{smallmatrix}
\right)  $, $h\in D\left(  \Gamma\right)  $, $h\left(  z\right)  =\sum
_{n=0}^{\infty}c_{n}z^{n}$, we get%
\begin{multline}
\left\langle k,Jk\right\rangle =2\operatorname{Re}\left\langle h,\Gamma
h\right\rangle =2\operatorname{Re}\left(  \sum_{n=0}^{\infty}\sum
_{m=0}^{\infty}\bar{c}_{n}\gamma_{n+m}c_{m}\right) \label{eqHan.17}\\
=2\sum_{n=0}^{\infty}\sum_{m=0}^{\infty}\bar{c}_{n}\operatorname{Re}\left(
\gamma_{n+m}\right)  c_{m}\geq0.
\end{multline}
But this means that the Hamburger moment problem is solvable for the sequence
$\left(  \operatorname{Re}\left(  \gamma_{n}\right)  \right)  _{n=0}^{\infty}%
$. If the solution is represented as in (\ref{eqHan.16}), then it follows that
$S_{\mu}\left(  U\right)  $ is represented as multiplication by $x$ on
$L^{2}\left(  \mathbb{R},\mu\right)  $, and we saw (using Theorem
\ref{ThmRef.1}(\ref{eqThmRef.1(5)})) that this forces $\mu$ to be supported in
the interval $I=\left[  -1,1\right]  $. Since $\gamma\in\ell^{2}$, it is known
from the theory of moments that $\mu$ is unique from $\Gamma_{\gamma}$. We
include the argument for why $S_{\mu}\left(  U\right)  $ is indeed
multiplication by $x$ on $L^{2}\left(  \mathbb{R},\mu\right)  $. Returning to
(\ref{eqHan.17}), we note that $S\left(  U\right)  $ is determined from the
identity%
\[
\left\langle k,JUk\right\rangle =\left\langle k,S\left(  U\right)
k\right\rangle _{J}%
\]
for $k=\left(
\begin{smallmatrix}
h\\
\Gamma h
\end{smallmatrix}
\right)  $, $h\in D\left(  \Gamma\right)  $; and we have:%
\[
JUk=%
\begin{pmatrix}
0 & I\\
I & 0
\end{pmatrix}%
\begin{pmatrix}
A & 0\\
0 & A^{\ast}%
\end{pmatrix}%
\begin{pmatrix}
h\\
\Gamma h
\end{pmatrix}
=%
\begin{pmatrix}
A^{\ast}\Gamma h\\
Ah
\end{pmatrix}
=%
\begin{pmatrix}
\Gamma Ah\\
Ah
\end{pmatrix}
.
\]
Consider finite sums $h_{1}\left(  z\right)  =\sum_{n}a_{n}z^{n}$ and
$h_{2}\left(  z\right)  =\sum_{n}b_{n}z^{n}$, and the corresponding
restrictions to $z=x\in\mathbb{R}$. Using $k_{1}=\left(
\begin{smallmatrix}
h_{1}\\
\Gamma h_{1}%
\end{smallmatrix}
\right)  $ and $k_{2}=\left(
\begin{smallmatrix}
h_{2}\\
\Gamma h_{2}%
\end{smallmatrix}
\right)  $, we get%
\begin{multline*}
\left\langle k_{1},S\left(  U\right)  k_{2}\right\rangle _{J}=\left\langle
h_{1},\Gamma Ah_{2}\right\rangle +\left\langle \Gamma h_{1},Ah_{2}%
\right\rangle =2\sum_{n}\sum_{m}\bar{a}_{n}\operatorname{Re}\left(
\gamma_{n+m}\right)  b_{m-1}\\
=\sum_{n}\sum_{m}\bar{a}_{n}\int_{\mathbb{R}}x^{n+m}\,d\mu\left(  x\right)
\,b_{m-1}=\int_{\mathbb{R}}\overline{h_{1}\left(  x\right)  }xh_{2}\left(
x\right)  \,d\mu\left(  x\right)  .
\end{multline*}
This concludes the proof of existence.\renewcommand{\qed}{}
\end{proof}

\begin{proof}
[Proof of uniqueness in Theorem \textup{\ref{ThmHan.3}}]Let $\mu$ be a finite
positive Borel measure on $\mathbb{R}$ which is supported in $\left[
-1,1\right]  $, and assume that $n\mapsto\int_{-1}^{1}x^{n}\,d\mu\left(
x\right)  $ is in $\ell^{2}$. We wish to reconstruct $\mathcal{K}=G\left(
\Gamma\right)  $ such that $\Gamma$ is a closed dissipative operator with
dense domain in $H^{2}$. Note that if $\Gamma$ has been found, then
\begin{equation}
\left\|
\begin{pmatrix}
h\\
\Gamma h
\end{pmatrix}
\right\|  _{J}^{2}=\left\langle h,\Gamma h\right\rangle +\left\langle \Gamma
h,h\right\rangle =\left\langle h,\left(  \Gamma+\Gamma^{\ast}\right)
h\right\rangle . \label{eqHan.a}%
\end{equation}
It follows that if $\Gamma\sim\left(  \gamma\right)  $ for some $\gamma\in
\ell^{2}$, then $\left\|  \left(
\begin{smallmatrix}
h\\
\Gamma h
\end{smallmatrix}
\right)  \right\|  _{J}$ and therefore the corresponding norm-completion
$\mathcal{H}_{J}\left(  G\left(  \Gamma\right)  \right)  $ only depends on the
sequence $\left(  \operatorname{Re}\gamma_{n}\right)  $, i.e., from
(\ref{eqHan.a}), $\Gamma+\Gamma^{\ast}\sim\left(  2\operatorname{Re}\gamma
_{n}\right)  $. Equivalently, we may assume without loss of generality that
the sequence $\left(  \gamma_{n}\right)  $ is real-valued. Now set
\begin{equation}
\gamma_{n}:=\frac{1}{2}\int_{-1}^{1}x^{n}\,d\mu\left(  x\right)  ,
\label{eqHan.b}%
\end{equation}
and let $\Gamma$ be the corresponding positive Hankel operator. For domain
$D\left(  \Gamma\right)  $, take the functions $h\in H^{2}$ which derive from
corresponding $\phi\in L_{\mu}^{2}\left(  \left[  -1,1\right]  \right)  $ as%
\begin{equation}
h\left(  z\right)  =\int_{-1}^{1}\left(  1-xz\right)  ^{-1}\phi\left(
x\right)  \,d\mu\left(  x\right)  . \label{eqHan.c}%
\end{equation}
Recall $A$ is the unilateral shift, and therefore%
\[
A^{\ast\,n}\gamma=\left(  \gamma_{n},\gamma_{n+1},\dots\right)  ,
\]
or, in function form,%
\begin{equation}
\left(  A^{\ast\,n}\gamma\right)  \left(  z\right)  =\gamma_{n}+\gamma
_{n+1}z+\gamma_{n+2}z^{2}+\cdots. \label{eqHan.d}%
\end{equation}
We then set
\begin{equation}
\left(  \Gamma h\right)  \left(  z\right)  =\sum_{n=0}^{\infty}\left(
A^{\ast\,n}\gamma\right)  \left(  z\right)  \int_{-1}^{1}x^{n}\phi\left(
x\right)  \,d\mu\left(  x\right)  \label{eqHan.e}%
\end{equation}
and note that $\Gamma$ is a Hankel operator, which is closed with dense domain
$D\left(  \Gamma\right)  \subset H^{2}$ and given by (\ref{eqHan.c}).
Moreover, $\mathcal{K}=G\left(  \Gamma\right)  $ has the desired properties,
with%
\begin{equation}
W_{\mu}%
\begin{pmatrix}
h\\
\Gamma h
\end{pmatrix}
\left(  x\right)  =h\left(  x\right)  \text{\qquad for }h\in D\left(
\Gamma\right)  \subset H^{2}, \label{eqHan.f}%
\end{equation}
and restricting $h$ to $\left(  -1,1\right)  \subset D$. Moreover, for
$\phi\in L_{\mu}^{2}\left(  \left[  -1,1\right]  \right)  $,
\begin{equation}
\left(  W_{\mu}^{\ast}\phi\right)  \left(  z\right)  =\int_{-1}^{1}\left(
1-xz\right)  ^{-1}\phi\left(  x\right)  \,d\mu\left(  x\right)
\label{eqHan.g}%
\end{equation}
is the function $h\left(  z\right)  $ given in (\ref{eqHan.c}) above.
\end{proof}

\begin{remark}
\label{RemHan154boundedness}\textup{(Boundedness)} The conditions
\textup{(\ref{LemHan154(1)(1)})--(\ref{LemHan154(1)(4)})} of Lemma
\textup{\ref{LemHan154}} are satisfied if $\gamma_{n}=\mathcal{O}\left(
\frac{1}{n}\right)  $, but, of course, for many examples which are not
$\mathcal{O}\left(  \frac{1}{n}\right)  $ as well. It is known in fact that
the Hankel operator $\Gamma_{\gamma}$ is bounded if and only if $\gamma
_{n}=\mathcal{O}\left(  \frac{1}{n}\right)  $. A theorem of Widom \cite{Wid66}
shows further that boundedness of the Hankel operator $\Gamma_{\gamma}$
\textup{(}from $\gamma_{n}=\int_{-1}^{1}x^{n}\,d\mu\left(  x\right)  $ with
$\mu$ a positive Borel measure\/\textup{)} holds if and only if $\mu$ is a
Carleson measure. \textup{(}A positive Borel measure $\mu$ on $I=\left[
-1,1\right]  $ is said to be a Carleson measure \cite{Car62} if and only if
$\mu\left(  I\setminus\left(  -x,x\right)  \right)  =\mathcal{O}\left(
1-x\right)  $ for $0<x<1$.\textup{)} It follows in particular that condition
\textup{(\ref{eqHan154.1})} in the Lemma \textup{\ref{LemHan154}} is satisfied
whenever $\Gamma_{\gamma}$ is assumed bounded; and further that
\textup{(\ref{eqHan154.1})} is more restrictive than requiring that $\left(
\gamma_{n}\right)  \in\ell^{2}$ where $\left(  \gamma_{n}\right)
_{n=0}^{\infty}$ denotes the moment sequence of $\mu$.
\end{remark}

\begin{remark}
\label{RemHan}The moment problem \textup{(\ref{eqHan.16})} \emph{with} the
finite support constraint seems to have been first studied in Devinatz
\cite[Lemma 1, p.~64]{Dev53}.
\end{remark}

\begin{corollary}
\label{CorHan.5}Let $\mathcal{K}=G\left(  \Gamma_{\gamma}\right)  $ be a
subspace of $H^{2}\oplus H^{2}$ satisfying the conditions in Theorem
\textup{\ref{ThmHan.3}.} Let $\mu$ be the measure on $\left[  -1,1\right]  $
given by%
\[
2\operatorname{Re}\gamma_{n}=\int_{-1}^{1}x^{n}\,d\mu\left(  x\right)
\qquad\left(  \in\ell^{2}\right)  ,
\]
and let%
\[
W_{\mu}\colon\mathcal{K}\longrightarrow L^{2}\left(  \left[  -1,1\right]
,\mu\right)
\]
be the contractive operator which intertwines $A\oplus A^{\ast}$ with
multiplication by $x$ on $L^{2}\left(  \left[  -1,1\right]  ,\mu\right)  $,
see Theorem \textup{\ref{ThmRef.1}.} Then%
\[
\ker\left(  W_{\mu}\right)  =\left\{  0\right\}
\]
if and only if $\operatorname*{supp}\left(  \mu\right)  $ has points of
accumulation in $\left(  -1,1\right)  $. \textup{(}So in particular, we can
have $\ker\left(  W_{\mu}\right)  =\left\{  0\right\}  $ both for measures
$\mu$ which are absolutely continuous relative to Lebesgue measure on $\left[
-1,1\right]  $, as well as for singular measures.\textup{)}
\end{corollary}

\begin{proof}
It follows from Theorem \ref{ThmHan.3} that
\begin{equation}
\left\|  W_{\mu}%
\begin{pmatrix}
h\\
\Gamma h
\end{pmatrix}
\right\|  ^{2}=\int_{-1}^{1}\left|  h\left(  x\right)  \right|  ^{2}%
\,d\mu\left(  x\right)  \text{,\qquad for }h\in H^{2}. \label{eqHan.18}%
\end{equation}
So for some $\left(  c_{n}\right)  \in\ell_{+}^{2}$, $h\left(  z\right)
=\sum_{n=0}^{\infty}c_{n}z^{n}$, and we may view $h\left(  x\right)  $ as the
restriction to $\left(  -1,1\right)  $ of the corresponding function $h\left(
z\right)  $ defined and analytic in $D=\left\{  z\in\mathbb{C}\mathrel
{;}\left|  z\right|  <1\right\}  $. If $\operatorname*{supp}\left(
\mu\right)  $ has accumulation points in $\left(  -1,1\right)  $, and $W_{\mu
}h=0$, then by (\ref{eqHan.18}), $h$ vanishes on a subset of
$\operatorname*{supp}\left(  \mu\right)  $ of full measure. This subset must
also have accumulation points, and since $h$ is analytic in $D$, it must
vanish identically.

To prove the converse, suppose $\operatorname*{supp}\left(  \mu\right)  $
contains only isolated points. Then $\mu$ must have the form
\[
\mu=\sum_{n}p_{n}\delta_{x_{n}},
\]
where $p_{n}>0$, and $\sum p_{n}<\infty$ and $\sum p_{n}\left(  1-x_{n}%
^{2}\right)  ^{-\frac{1}{2}}<\infty$. Recall $\mu$ is finite, and supported in
$\left[  -1,1\right]  $. Then pick $h\in H^{2}$, $\left\|  h\right\|  _{H^{2}%
}\neq0$, such that $h\left(  x_{n}\right)  =0$, for example $h\left(
z\right)  =\left(  \prod_{n}\frac{x_{n}-z}{1-x_{n}z}\right)  z\left(
1-z^{2}\right)  $. Then $h\in\ker\left(  W_{\mu}\right)  $.
\end{proof}

\begin{acknowledgements}
The problems addressed in the present paper grew out of earlier joint work
with G. \'{O}lafsson \cite{JoOl98,JoOl99} as well as earlier work by the
present author. We are very grateful to G. \'{O}lafsson for the benefit of
ongoing discussions. We are also very grateful to Brian Treadway for excellent
manuscript production.
\end{acknowledgements}

\bibliographystyle{bftalpha}
\bibliography{jorgen}
\end{document}